\documentclass[12pt]{article}
\usepackage{latexsym, graphicx, amscd}
\DeclareGraphicsExtensions{.ps}

    \title{{\bf  A functional-analytic
theory of vertex  (operator)
algebras, I}}
    \author{Yi-Zhi Huang}
    \date{}
    
\begin{document}
    \bibliographystyle{alpha}
    \maketitle

    \input amssym.def
    \input amssym
    \newtheorem{rema}{Remark}[section]
    \newtheorem{propo}[rema]{Proposition}
    \newtheorem{theo}[rema]{Theorem} 

   \newtheorem{defi}[rema]{Definition}
    \newtheorem{lemma}[rema]{Lemma}
    \newtheorem{corol}[rema]{Corollary}
     \newtheorem{exam}[rema]{Example}
\newcommand{\binom}[2]{{{#1}\choose {#2}}}
	\newcommand{\nno}{\nonumber}
	\newcommand{\lbar}{\bigg\vert}
\newcommand{\mbar}{\mbox{\large $\vert$}}
	\newcommand{\p}{\partial}
	\newcommand{\dps}{\displaystyle}
	\newcommand{\bra}{\langle}
	\newcommand{\ket}{\rangle}
 \newcommand{\res}{\mbox{\rm Res}}
\renewcommand{\hom}{\mbox{\rm Hom}}
 \newcommand{\pf}{{\it Proof.}\hspace{2ex}}
 \newcommand{\epf}{\hspace{2em}$\Box$}
 \newcommand{\epfv}{\hspace{2em}$\Box$\vspace{1em}}
 \newcommand{\epfe}{\hspace{2em}\Box}
\newcommand{\nord}{\mbox{\scriptsize ${\circ\atop\circ}$}}
\newcommand{\wt}{\mbox{\rm wt}\ }
\newcommand{\swt}{\mbox{\rm {\scriptsize wt}}\ }
\newcommand{\clr}{\mbox{\rm clr}\ }

\begin{abstract}
This paper is the first in a series of papers
developing a functional-analytic theory of vertex
(operator) algebras and their representations. 
For an arbitrary ${\Bbb Z}$-graded finitely-generated
vertex algebra $(V, Y, {\bf 1})$ satisfying the standard 
grading-restriction axioms, 
a locally convex topological completion $H$ of 
$V$ is constructed.  By the geometric interpretation of 
vertex (operator) algebras,  there is a canonical linear
map from
$V\otimes V$ to $\overline{V}$ (the algebraic completion of $V$)
realizing linearly
the conformal equivalence class 
of a genus-zero Riemann 
surface with 
analytically parametrized boundary obtained by deleting
two ordered disjoint disks
from the unit disk and by giving the obvious parametrizations to the 
boundary components.
We extend such a linear map to a linear map from $H\tilde{\otimes} H$
($\tilde{\otimes}$ being the completed
tensor product) to $H$,
and prove the continuity of the extension.
For any finitely-generated ${\Bbb C}$-graded 
$V$-module $(W, Y_{W})$ 
satisfying the standard grading-restriction axioms,
the same method also gives a topological completion $H^{W}$
of $W$ 
and gives the continuous extensions from $H\tilde{\otimes}H^{W}$ 
to $H^{W}$ of the linear maps from
$V\otimes 
W$ to $\overline{W}$ 
realizing linearly the above conformal equivalence 
classes of the genus-zero Riemann 
surfaces with 
analytically parametrized boundaries.
\end{abstract}

\renewcommand{\theequation}{\thesection.\arabic{equation}}
\renewcommand{\therema}{\thesection.\arabic{rema}}
\setcounter{equation}{0}
\setcounter{rema}{0}
\setcounter{section}{-1}

\section{Introduction}

We begin a systematic study of  the functional-analytic structure of vertex
(operator) algebras and their representations in this paper. 

Vertex (operator) algebras were introduced rigorously in mathematics
by Borcherds and by Frenkel, Lepowsky and Meurman (see \cite{B},
\cite{FLM} and \cite{FHL}). Incorporating modules and intertwining
operators, the author introduced intertwining operator algebras in
\cite{H5}.  The original definition of vertex (operator)
algebra is purely algebraic,
but a geometric interpretation of vertex operators motivated by the
path integral picture in string theory was soon observed
mathematically by Frenkel \cite{F}.  In \cite{H1}--\cite{H7}, the
author developed a geometric theory of vertex operator algebras and
intertwining operator algebras; the hard parts deal with the Virasoro
algebra, the central charge, the interaction between the Virasoro
algebra and vertex operators (or intertwining operators), and the
monodromies. But in this geometric theory, the linear maps associated to
elements of certain moduli spaces of punctured spheres with local
coordinates (or to elements of the vector bundles forming partial operads
called ``genus-zero 
modular functors'' over these moduli spaces) are from
the tensor powers of a vertex operator algebra (or an intertwining
operator algebra) to the algebraic completion of the algebra. 

To be more specific, let $(V, Y, {\bf 1})$ be a ${\Bbb Z}$-graded
vertex algebra. 
For any nonzero complex number $z$, let $P(z)$ be the 
the conformal equivalence class 
of the sphere ${\Bbb C}\cup \{\infty\}$ with the negatively oriented 
puncture $\infty$, the ordered positively oriented $z$ and $0$, and 
with the standard local coordinates. Then associated to $P(z)$
is the linear map
$Y(\cdot, z)\cdot: V\otimes V\to \overline{V}$ given by the vertex operator
map. Note that the image is  in 
the algebraic completion $\overline{V}$ of $V$, not $V$ itself.

One of the main goals of the theory of vertex operator algebras is to
construct conformal field theories in the sense of Segal \cite{S1}
\cite{S2} from vertex operator algebras and their representations, and
to study geometry using the theory of vertex operator algebras.  It is
therefore necessary to construct locally convex topological
completions of the underlying vector spaces of vertex (operator)
algebras and their modules, such that associated to any Riemann
surface with boundary, one can construct continuous linear maps
between the tensor powers of these completions. In particular, it is necessary 
to construct a locally convex completion $H$ of the underlying vector space
$V$ of a vertex algebra $(V, Y, {\bf 1})$ such that associated to the 
conformal equivalence class of 
a disk with two smaller ordered disks deleted and with 
the obvious boundary parametrization, there is a canonical continuous map 
from $H\tilde{\otimes}H$ ($\tilde{\otimes}$ being the completed
tensor product) to $H$. 
Though in some algebraic applications of conformal field theories,
only the algebraic structure of conformal field theories 
is needed, in many geometric applications, 
it is  necessary to have a complete
locally convex topological space and continuous linear maps, associated
to Riemann surfaces with boundaries, between tensor powers of the
space.

In the present paper (Part I), we 
construct a locally convex completion $H$ of an arbitrary finitely-generated
${\Bbb Z}$-graded vertex algebra $(V, Y, {\bf 1})$
satisfying the standard grading-restriction
axioms.
The completion $H$ is the strict inductive limit of a sequence of complete
locally convex spaces constructed from the correlation functions 
of the generators. The strong topology, rather than the weak-$*$
topology, on the topological dual spaces of certain function spaces 
is needed in the construction.
For any positive number
$r_{1}, r_{2}$ and nonzero complex number $z$ satisfying
$r_{2}+2r_{1}<1$ and $r_{2}<|z|<1$, there is a unique genus-zero Riemann 
surface with 
analytically parametrized boundary given by deleting
two ordered disjoint disks, the first centered at $z$ with 
radius $r_{1}$ and the second centered at $0$ with radius $r_{2}$, 
from the unit disk and by giving the obvious parametrizations to the 
boundary components.
Associated to this genus-zero Riemann surface with 
analytically parametrized boundary is a linear map 
$Y(r_{1}^{L(0)}\cdot, z)r_{2}^{L(0)}: V\otimes V\to \overline{V}$.
We extend this linear map to a linear map from $H\tilde{\otimes} H$
($\tilde{\otimes}$ being the completed
tensor product) to $H$,
and prove the continuity of the extension.
It is clear that $H$ is linearly isomorphic
to a subspace of $\overline{V}$ containing both $V$ and the image of 
$Y(r_{1}^{L(0)}\cdot, z)r_{2}^{L(0)}$.

For any finitely-generated ${\Bbb C}$-graded 
module $(W, Y_{W})$ satisfying the standard grading-restriction axioms 
for such a finitely-generated 
vertex algebra,
we also construct a locally convex completion $H^{W}$
of $W$, and construct continuous extensions from $H\tilde{\otimes}H^{W}$ 
to $H^{W}$ of the linear map 
$Y_{W}(r_{1}^{L(0)}\cdot, z)r_{2}^{L(0)}: V\otimes 
W\to \overline{W}$.

Our method depends only on the
axiomatic properties or ``world-sheet geometry'' (mainly the duality
properties) of the vertex  algebra. Since our construction
does not use any additional structure, we expect that the locally
convex 
completions constructed in the present paper will be useful in
solving purely algebraic problems in the representation theory of 
finitely-generated 
vertex (operator) algebras.

This paper is organized as follows: In Section 2, we construct a
locally convex completion $H$ of a finitely-generated ${\Bbb Z}$-graded 
vertex algebra $(V, Y, {\bf 1})$ satisfying the standard
grading-restriction axioms.
In Section 3, we extend $Y(r_{1}^{L(0)}\cdot, z)r_{2}^{L(0)}: V\otimes V\to 
\overline{V}$ to continuous linear maps from $H\tilde{\otimes} H$
to $H$. In Section 4, we
present the corresponding results for modules. 

In this paper, we assume that the reader is 
familiar with the basic definitions and 
results in the theory of vertex algebras. The material in 
\cite{B}, \cite{FLM} and \cite{FHL} should be enough. 
We also assume that the reader is familiar with the basic definitions, 
constructions and results in the theory of locally convex topological vector
spaces. The reader can find this material in, for example, 
\cite{K1} and \cite{K2}.

We shall denote the set of integers, the set of real numbers and the set
of complex numbers by the usual notations ${\Bbb Z}$, ${\Bbb R}$ and
${\Bbb C}$, respectively.  We shall use $i, j, k, l, m, n, p, q$
to denote integers. In particular, when we write, say, $k>0$ (or $k\ge
0$), 
we mean 
that $k$ is a positive integer (or a nonnegative integer). For a
graded vector space $V$, we use $V'$, $V^{*}$ and $\overline{V}$ to
denote the graded dual space, the dual space and the algebraic
completion of $V$, respectively. For a topological vector space $E$,
we use $E^{*}$ to denote the topological dual space of $E$. The symbol
$\otimes$ always denotes the vector space tensor product. The bifunctor
given by
completing the vector space tensor product of two topological
vector spaces with the tensor product topology is denoted by
$\tilde{\otimes}$.

\paragraph{Acknowledgment} This research is supported in part by 
by NSF grant DMS-9622961.

\renewcommand{\theequation}{\thesection.\arabic{equation}}
\renewcommand{\therema}{\thesection.\arabic{rema}}
\setcounter{equation}{0}
\setcounter{rema}{0}

\section{A locally convex completion of
a finitely-gen\-erated vertex  algebra}

In this section, we use ``correlation functions'' to
construct the locally convex completion of a finitely-generated vertex
algebra. 

For any $k\ge 0$, let $R_{k}$ be the space of rational functions 
in the complex variables $z_{1}, \dots, z_{k}$ with the only
possible poles 
$z_{i}=z_{j}$ for $i\ne j$ and
$z_{i}=0, \infty$ ($i, j=1, \dots, k$). 
Let 
$$M^{k}=\{(z_{1}, \dots, z_{k})
\in {\Bbb C}^{k}\;|\; 
z_{i}\ne z_{j}\;\;\mbox{for}\;\; i\ne j; z_{i}\ne 0\;(i, j=1, \dots, k)\}$$
and let 
$\{K^{(k)}_{n}\}_{n>0}$, be a sequence of compact subsets of 
$M^{k}$ satisfying $$K^{(k)}_{n}\subset K^{(k)}_{n+1}, \;\;n>0,$$
and 
$$M^{k}=\bigcup_{n>0}K^{(k)}_{n}.$$
For any $n>0$, we define a map
$\|\cdot\|_{R_{k}, n}: R_{k}\to [0, \infty)$
by 
$$\|f\|_{R_{k}, n}=\max_{(z_{1}, \dots, z_{k})\in K^{(k)}_{n}}
|f(z_{1}, \dots, z_{k})|$$
for $f\in R_{k}$. Then it is clear that $\|\cdot\|_{R_{k}, n}$ is a 
norm on $R_{k}$. Using this sequence of norms, we obtain a locally convex
topology on $R_{k}$.  Note that 
a sequence in $R_{k}$ is convergent if and only if this sequence of
functions is uniformly convergent on any compact subset of $M^{k}$.
Clearly, this topology is independent of the 
choice of the sequence $\{K^{(k)}_{n}\}_{n>0}$.

Let $V$ be a ${\Bbb Z}$-graded vertex algebra satisfying 
the standard grading-restriction axioms, that is, 
$$V=\coprod_{n\in {\Bbb Z}}V_{(n)},$$
$$\dim V_{(n)}<\infty$$
for $n\in {\Bbb Z}$
and 
$$V_{(n)}=0$$
for $n$ sufficiently small. 
By the duality properties 
of $V$,
for any $v'\in
V'$, any $u_{1}, \dots, u_{k}, v\in V$, 
$$\langle v', Y(u_{1},
z_{1})\cdots Y(u_{k}, z_{k})v \rangle$$ 
is absolutely convergent in the region
$|z_{1}|>\cdots>|z_{k}|>0$ 
and can be analytically extended to an element
$$
R(\langle v', Y(u_{1},
z_{1})\cdots Y(u_{k}, z_{k})v \rangle)
$$
of $R_{k}$.

For any $u_{1}, \dots, u_{k}, v\in V$ and any $(z_{1}, \dots, 
z_{k})\in M^{k}$, we have an element 
$$Q(u_{1}, \dots, u_{k}, v; z_{1}, \dots, z_{k})\in \overline{V}$$
defined by
$$\langle v',
Q(z_{1}, \dots, z_{k}, v; z_{1}, \dots, z_{k})\rangle
=R(\langle v', Y(u_{1},
z_{1})\cdots Y(u_{k}, z_{k})v \rangle)$$
for $v'\in V'$. We denote the projections from $V$ to $V_{(n)}$, $n\in
+{\Bbb Z}$, by $P_{n}$.
Let $\tilde{G}$ be the subspace of $V^{*}$ 
consisting of linear functionals $\lambda$ on $V$
such that for any $k\ge 0$,
$u_{1}, 
\dots, u_{k}, v\in V$,
\begin{equation}\label{conv-fk}
\sum_{n\in {\Bbb Z}}\lambda(P_{n}(Q(u_{1}, \dots, u_{k}, v; z_{1}, \dots,
z_{k})))
\end{equation}
is absolutely convergent for any $z_{1}, \dots, z_{k}$
in the region 
$$M^{k}_{<1}=\{(z_{1}, \dots, z_{k})\in M^{k}\;|\;
|z_{1}|, \dots, |z_{k}|<1\}.$$
The dual pair $(V^{*}, V)$  of vector spaces
gives $V^{*}$ a locally convex topology.
With the topology induced from the one on $V^{*}$,
$\tilde{G}$ is also a 
locally convex  space. Note that $V'$ is a subspace of $\tilde{G}$.

For any $k\ge 0$, 
$\lambda\in \tilde{G}$, $u_{1}, 
\dots, u_{k}, v\in V$, both
(\ref{conv-fk}) and
\begin{eqnarray*}
\lefteqn{\sum_{n\in {\Bbb Z}}\frac{\p}{\p z_{1}}
\lambda(P_{n}(Q(u_{1}, \dots, u_{k}, v; z_{1}, \dots,
z_{k})))}\nno\\
&&=\sum_{n\in {\Bbb Z}}\lambda(P_{n}(Q(L(-1)u_{1}, 
\dots, u_{k}, v; z_{1}, \dots,
z_{k})))
\end{eqnarray*}
are absolutely convergent in the region $M^{k}_{<1}$.
Thus (\ref{conv-fk}) is analytic in $z_{1}$ 
when $(z_{1}, \dots, z_{k})$ is in $M^{k}_{<1}$.
Similarly, (\ref{conv-fk}) is analytic in $z_{i}$ for 
$i=2, \dots, k$ when $(z_{1}, \dots, z_{k})$
is in $M^{k}_{<1}$. So (\ref{conv-fk})
defines an analytic function on $M^{k}_{<1}$ and we denote it by
$$g_{k}(\lambda \otimes u_{1}\otimes \cdots\otimes u_{k}
\otimes v)$$
since (\ref{conv-fk}) is multilinear in $\lambda$, $u_{1}, \dots, u_{k}$
and $v$. 
These functions  span a vector space $F_{k}$ of analytic functions on 
$M^{k}_{<1}$. So we obtain a linear map 
$$g_{k}: \tilde{G}\otimes V^{\otimes (k+1)}\to F_{k}.$$
Fix a sequence $\{J^{(k)}_{n}\}_{n>0}$
of compact subsets 
of $M^{k}_{<1}$ such that 
$$J^{(k)}_{n}\subset J^{(k)}_{n+1},\;\;
n>0,$$
and 
$$\bigcup_{n>0}J^{(k)}_{n}=M^{k}_{<1}.$$
As in the case of
$R_{k}$, using these compact subsets, we define a sequence of 
norms $\|\cdot\|_{F_{k}, n}$ on $F_{k}$, and these norms 
give a locally convex topology on $F_{k}$. 

There is also an embedding $\iota_{F_{k}}$ from
$F_{k}$ to $F_{k+1}$ defined as follows:  We use $(z_{0}, \dots,
z_{k})$ instead of $(z_{1}, \dots, z_{k+1})$ to denote the elements of
$M^{k+1}_{<1}$. 
For $\lambda\in \tilde{G}$, $u_{1}, \dots, u_{k}, v\in V$, since 
$$Y({\bf 1}, z)=1$$
for any nonzero complex number $z$,
$$g_{k+1}(\lambda \otimes  {\bf
1}\otimes u_{1}\otimes \cdots \otimes u_{k}\otimes 
v)$$ as a function of $z_{0}, \dots, z_{k}$
is in fact independent of $z_{0}$, and is equal to
$$g_{k}(\lambda \otimes u_{1}\otimes \cdots \otimes u_{k}\otimes 
v)$$ 
as a function in $z_{1}, \dots, z_{k}$. Thus we obtain a
well-defined linear
map 
$$\iota_{F_{k}}: F_{k}\to F_{k+1}$$
such that 
$$\iota_{F_{k}}\circ g_{k}=g_{k+1}\circ \phi_{k}$$
where 
$$\phi_{k}: \tilde{G}\otimes V^{\otimes (k+1)}  \to \tilde{G} \otimes 
V^{\otimes (k+2)}$$
is defined by
$$\phi_{k}(\lambda \otimes  u_{1}\otimes \cdots \otimes u_{k}\otimes v)
=\lambda\otimes {\bf 1}\otimes 
  u_{1}\otimes \cdots \otimes u_{k}\otimes v$$
for $\lambda\in \tilde{G}$, $u_{1}, \dots, u_{k},
v\in V$. It is clear
that $\iota_{F_{k}}$ is injective. Thus we can regard $F_{k}$ as a
subspace of $F_{k+1}$. Moreover, we have:

\begin{propo}\label{k:k+1}
For any $k\ge 0$,   $\iota_{F_{k}}$ as a map from $F_{k}$ to 
$\iota_{F_{k}}(F_{k})$
is continuous and open. In other words, the topology on $F_{k}$
is induced from that on $F_{k+1}$.
\end{propo}
\pf
We consider the two topologies on $F_{k}$, one is the topology
defined above for $F_{k}$ 
and the other induced from the topology on $F_{k+1}$.
We need only prove that for any $n>0$, (i) the norm
$\|\cdot\|_{F_{k}, n}$ is continuous in the topology induced from the one
on $F_{k+1}$, and (ii) the restriction of the norm
$\|\cdot\|_{F_{k+1}, n}$ to $F_{k}$ is continuous in the topology on
$F_{k}$. 

Let 
$\{f_{\alpha}\}_{\alpha\in A}$ be a net in $F_{k}$
convergent in the topology induced from the one on $F_{k+1}$.  Then
$\{f_{\alpha}\}_{\alpha\in A}$, when viewed as a net of
functions in $z_{0}, z_{1}, \dots, z_{k}$, 
is convergent uniformly
on any compact subset of $M^{k+1}$. Since $f_{\alpha}$, $\alpha\in A$, 
are independent of $z_{0}$, $\{f_{\alpha}\}_{\alpha\in A}$ 
is in fact convergent uniformly on any compact subset of
$M^{k}$, proving (i). Now let $\{f_{\alpha}\}_{\alpha\in A}$ be a
net in $F_{k}$ convergent in the topology on $F_{k}$. Then
$\{f_{\alpha}\}_{\alpha\in A}$ is convergent uniformly on any
compact subset of $M^{k}$. If we view $f_{\alpha}$, $\alpha\in A$,
as functions on ${\Bbb C}\times M^{k}$, then the net
$\{f_{\alpha}\}_{\alpha\in A}$ is convergent uniformly on any
subset of $M^{k+1}$ of the form ${\Bbb C}\times K$ where $K$ is a
compact subset of $M^{k}$. But any compact subset of $M^{k+1}$ is 
contained in a
subset of the form ${\Bbb C}\times K$. So $\{f_{\alpha}\}_{\alpha\in A}$
is convergent uniformly on any compact subset of $M^{k+1}$,
proving (ii). 
\epfv

We equip the topological dual space $F_{k}^{*}$, $k\ge 0$, of $F_{k}$
with the strong topology, that is, the topology of uniform convergence 
on all the weakly bounded subsets of $F_{k}$ (see page 256 of \cite{K1}).
Then $F_{k}^{*}$ is a locally convex space. In fact, $F_{k}^{*}$ is a
(DF)-space (see page 396 of \cite{K1}).

For any $\lambda\in \tilde{G}$ and $u\in V$, let
$Y_{-1}(u) \lambda\in \tilde{G}$ be defined by 
$$(Y_{-1}(u)\lambda)(v)=\lambda((\res_{x}x^{-1}Y(u, x))v)$$
for $v\in V$. (Note that $\res_{x}x^{-1}Y(u, x)$ is 
$u_{-1}$ in the usual notation. So we could denote $Y_{-1}(u)$ by $u_{-1}$.
Here we use the notation $Y_{-1}(u)$ to avoid the possible confusion
with the notations $u_{0}, u_{1}, \dots$ for elements in $V$ used later.)
For  $k\ge 0$, we define a linear map 
$$\gamma_{k}: F_{k+1}\to F_{k}$$
by
\begin{eqnarray*}
\lefteqn{\gamma_{k}(g_{k+1}(\lambda\otimes 
u_{0}\otimes u_{1}\otimes \cdots \otimes 
u_{k}\otimes v))}\nno\\
&&=g_{k}(Y_{-1}(u_{0}) \lambda)\otimes  u_{1}\otimes \cdots \otimes 
u_{k}\otimes v)
\end{eqnarray*}
for
$\lambda\in \tilde{G}$, $u_{0}, u_{1}, \dots, u_{l}, v\in V$.

\begin{propo}
The map $\gamma_{k}$ is continuous and satisfies 
\begin{equation}\label{halfinv}
\gamma_{k}
\circ \iota_{F_{k}}=I_{F_{k}}
\end{equation}
where $I_{F_{k}}$ is the identity map
on $F_{k}$. 
\end{propo}
\pf
We use $(z_{0}, \dots, z_{k})$ instead of $(z_{1}, \dots, z_{k+1})$ to
denote the elements of $M^{k+1}_{<1}$.
From the definition, we see that for any 
positive number $\epsilon<1$, when $|z_{1}|, \dots, |z_{k}|<\epsilon$,
\begin{eqnarray*}
\lefteqn{\gamma_{k}(g_{k+1}(\lambda\otimes
u_{0}\otimes u_{1}\otimes \cdots \otimes 
u_{l}\otimes v))}\nno\\
&&=
\frac{1}{2\pi \sqrt{-1}}\oint_{|z_{0}|=\epsilon}z_{0}^{-1}
g_{k+1}(\lambda\otimes u_{0}\otimes u_{1}\otimes \cdots \otimes 
u_{k}\otimes v)dz_{0}
\end{eqnarray*}
for $\lambda\in \tilde{G}$, $u_{0}, \dots, u_{k}, 
v\in V$. Thus for any $n>0$,
there exist $m_{n}>0$ and positive number $\epsilon_{n}$
such that 
\begin{eqnarray*}
\lefteqn{\|\gamma_{k}(g_{k+1}(\lambda\otimes
 u_{0}\otimes u_{1}\otimes \cdots \otimes 
u_{k}\otimes v))\|_{F_{k}, n}}\nno\\
&&=\max_{(z_{1}, \dots, z_{k})\in J_{n}^{(k)}}
|\gamma_{k}(g_{k+1}(\lambda\otimes u_{0}\otimes u_{1}\otimes \cdots \otimes 
u_{k}\otimes v))|\nno\\
&&=\max_{(z_{1}, \dots, z_{k})\in J_{n}^{(k)}}
\biggl|\frac{1}{2\pi \sqrt{-1}}\oint_{|z_{0}|=\epsilon_{n}}z_{0}^{-1}
g_{k+1}(\lambda\otimes u_{0}\otimes u_{1}\otimes \cdots \otimes 
u_{k}\otimes v)dz_{0}\biggr|\nno\\
&&\le \max_{(z_{1}, \dots, z_{k})\in J_{n}^{(k)}, |z_{0}|=\epsilon_{n}}
|g_{k+1}(\lambda\otimes u_{0}\otimes u_{1}\otimes \cdots \otimes 
u_{k}\otimes v)|\nno\\
&&\le \max_{(z_{0}, \dots, z_{k})\in J_{m_{n}}^{(k+1)}}
|g_{k+1}(\lambda\otimes u_{0}\otimes u_{1}\otimes \cdots \otimes 
u_{k}\otimes v)|\nno\\
&&=\|g_{k+1}(\lambda\otimes  u_{0}\otimes u_{1}\otimes \cdots \otimes 
u_{k}\otimes v)\|_{F_{k+1}, m_{n}}.
\end{eqnarray*}
This inequality implies that $\gamma_{k}$ is continuous.

For $\lambda\in \tilde{G}$, $u_{1}, \cdots,
u_{k}, v\in V$, by definition,
\begin{eqnarray*}
\lefteqn{g_{k+1}(\lambda\otimes  {\bf 1}\otimes u_{1}\otimes \cdots \otimes 
u_{k}\otimes v)}\nno\\
&&=\iota_{F_{k}}(g_{k}(\lambda
  \otimes u_{1}\otimes \cdots \otimes 
u_{k}\otimes v)).
\end{eqnarray*}
Since $Y_{-1}({\bf 1}) \lambda=\lambda$,
\begin{eqnarray*}
\lefteqn{\gamma_{k}(g_{k+1}(\lambda\otimes {\bf 1}
\otimes u_{1}\otimes \cdots \otimes 
u_{k}\otimes v))}\nno\\
&&=g_{k}((Y_{-1}({\bf 1}) \lambda)
  \otimes u_{1}\otimes \cdots \otimes 
u_{k}\otimes v)\nno\\
&&=g_{k}(\lambda
  \otimes u_{1}\otimes \cdots \otimes 
u_{k}\otimes v).
\end{eqnarray*}
 So we have
(\ref{halfinv}).
\epfv

\begin{corol}\label{k*:k+1*}
The  adjoint 
map $\gamma_{k}^{*}$  of $\gamma_{k}$ 
satisfies 
\begin{equation}\label{halfinv*}
\iota_{F_{k}}^{*}\circ \gamma_{k}^{*}
=I_{F_{k}^{*}}
\end{equation}
where  
$$\iota_{F_{k}}^{*}: 
F_{k+1}^{*}\to F_{k}^{*}$$
is the adjoint  of $\iota_{F_{k}}$ and 
$I_{F_{k}^{*}}$ is the identity on $F_{k}^{*}$. It  
is
injective and continuous. As a map from $F_{k}^{*}$ to
$\gamma_{k}^{*}(F_{k}^{*})$, it is also open.
In particular, if we identify $F_{k}^{*}$ 
with $\gamma_{k}^{*}(F_{k}^{*})$,
the topology on $F_{k}^{*}$ is 
induced from the one on $F_{k+1}^{*}$.
\end{corol}
\pf
The identity 
(\ref{halfinv*}) is an immediate consequence of the identity
(\ref{halfinv}). 
By this identity, we see that $\gamma_{k}^{*}$ is  injective. 
The continuity of $\gamma_{k}^{*}$
is a consequence of the continuity of
$\gamma_{k}$. To show that it is open as a map from $F_{k}^{*}$ to
$\gamma_{k}^{*}(F_{k}^{*})$, we need only show that 
its inverse  is continuous.
Let $\{\mu_{\alpha}\}_{\alpha\in A}$ 
be a net in $F_{k}^{*}$ such that 
$\{\gamma_{k}^{*}(\mu_{\alpha})\}_{\alpha\in A}$ 
is convergent in  $F_{k+1}^{*}$. 
Since $\iota_{F_{k}}^{*}$ is continuous,
$\{\iota_{F_{k}}^{*}(\gamma_{k}^{*}(\mu_{\alpha}))\}_{\alpha\in A}$ 
is convergent in $F_{k}^{*}$. By the identity (\ref{halfinv}),
$$\mu_{\alpha}=\iota_{F_{k}}^{*}(\gamma_{k}^{*}(\mu_{\alpha}))$$
for $\alpha\in A$. 
Thus $\{\mu_{\alpha}\}_{\alpha\in A}$  is convergent 
in $F_{k}^{*}$, proving that the  inverse of 
$\gamma_{k}^{*}$ viewed as a map from $F_{k}^{*}$ to
$\gamma_{k}^{*}(F_{k}^{*})$ is  continuous. 
\epfv

We use $\langle \cdot, \cdot\rangle$ to denote the pairing between
$\tilde{G}$ and the algebraic dual space $\tilde{G}^{*}$ of 
$\tilde{G}$. Since $V'\subset \tilde{G}$ and $V\subset \tilde{G}^{*}$,
this pairing is an extension of the pairing between $V'$ and $V$
denoted using the same symbol. 
With this pairing, $\tilde{G}$ and $\tilde{G}^{*}$ form a dual pair 
of vector spaces.
This dual pair of vector spaces 
gives a locally convex topology to $\tilde{G}^{*}$.
Since $V'\subset \tilde{G}$, 
the dual space $\tilde{G}^{*}$ can be viewed as a subspace of 
$(V')^{*}=\overline{V}$. 
We define 
$$e_{k}: V^{\otimes (k+1)}\otimes F_{k}^{*}\to
\tilde{G}^{*}\subset \overline{V}$$
by
$$
\langle \lambda, e_{k}(u_{1}\otimes \cdots \otimes u_{k}\otimes
v\otimes \mu)\rangle
=\mu(g_{k}(\lambda\otimes u_{1}\otimes \cdots \otimes u_{k}\otimes
v))
$$
for $\lambda\in \tilde{G}$, $u_{1}, \dots, u_{k}, v\in V$ and 
$\mu \in F_{k}^{*}$. 

We now assume that $V$ is finitely generated.
The generators of $V$ 
span a finite-dimensional subspace $X$ of $V$.
We assume that the vacuum vector 
${\bf 1}$ is in $X$.  Any norm on $X$ 
induces a Banach space structure on $X$. 
Since $X$ is finite-dimensional, all norms on $X$ are
equivalent, so that the topology induced by the norm  
is in fact independent of the norm. 

Since $X$ is a
finite-dimensional Banach space and $F_{k}^{*}$ is a 
locally convex space,
$X^{\otimes (k+1)}\otimes F^{*}_{k}$ is also a 
locally convex
space. We denote the image $e_{k}(X^{\otimes (k+1)}\otimes F_{k}^{*})$ of 
$X^{\otimes (k+1)}\otimes F_{k}^{*}\subset V^{\otimes (k+1)}\otimes 
F_{k}^{*}$
under $e_{k}$ by $G_{k}$.

\begin{propo}\label{inclusion}
For any $k\ge 0$, $G_{k}\subset G_{k+1}$.
\end{propo}
\pf 
By definition,
\begin{eqnarray*}
\lefteqn{\langle \lambda, e_{k}(u_{1}\otimes \cdots \otimes u_{k}\otimes
v\otimes \mu)\rangle}\nno\\
&&=\mu(g_{k}(\lambda\otimes u_{1}\otimes \cdots \otimes u_{k}\otimes
v)\nno\\
&&=\mu(\gamma_{k}(g_{k+1}(\lambda\otimes {\bf 1}\otimes 
u_{1}\otimes \cdots \otimes u_{k}\otimes
v))\nno\\
&&=(\gamma_{k}^{*}(\mu))(g_{k+1}(\lambda\otimes {\bf 1}\otimes 
u_{1}\otimes \cdots \otimes u_{k}\otimes
v))\nno\\
&&=\langle \lambda, e_{k}({\bf 1}\otimes 
u_{1}\otimes \cdots \otimes u_{k}\otimes
v\otimes \gamma_{k}^{*}(\mu))\rangle
\end{eqnarray*}
for $\lambda\in \tilde{G}$, $u_{1}, \dots, u_{k}, 
v\in V$ and $\mu \in F_{k}^{*}$.
Thus 
$$e_{k}(u_{1}\otimes \cdots \otimes u_{k}\otimes
v\otimes \mu)=e_{k}({\bf 1}\otimes 
u_{1}\otimes \cdots \otimes u_{l}\otimes
v\otimes \gamma_{k}^{*}(\mu))
\in G_{k+1}$$
for $\lambda\in \tilde{G}$, $u_{1}, \dots, u_{k}, 
v\in V$ and $\mu \in F_{k}^{*}$.
\epfv

\begin{propo}\label{elcont}
The linear map 
$$e_{k}|_{X^{\otimes (k+1)}\otimes F_{k}^{*}}: 
X^{\otimes (k+1)}\otimes F_{k}^{*}\to
\tilde{G}^{*}$$
is continuous. 
\end{propo}
\pf
By the definition of the locally
convex topology on
the tensor product $X^{\otimes (k+1)}\otimes
F_{k}^{*}$, we need only prove that $e_{k}|_{X^{\otimes (k+1)}\otimes
F_{k}^{*}}$ as a multilinear map from 
$X^{\otimes (k+1)}\times
F_{k}^{*}$ to $\tilde{G}^{*}$ is continuous. 

 Let 
 $\{({\cal X}_{\alpha}, 
\mu_{\alpha})\}_{\alpha\in A}$ be a net in $X^{\otimes (k+1)}\times
F_{k}^{*}$ convergent to $0$. Then the nets 
$\{{\cal X}_{\alpha}\}_{\alpha\in A}$
and $\{\mu_{\alpha}\}_{\alpha\in A}$ are convergent to $0$ in 
$X^{\otimes (k+1)}$ and $F_{k}^{*}$, respectively. 
For any fixed $\lambda\in \tilde{G}$, since $X^{\otimes (k+1)}$ 
is a finite-dimensional Banach space, the linear map from
$X^{\otimes (k+1)}$ to $F_{k}$ defined by 
$$u_{1}\otimes \cdots \otimes u_{k}\otimes v\mapsto 
g_{k}(\lambda\otimes u_{1}\otimes \cdots \otimes u_{k}\otimes v)$$
for $u_{1}, \dots, u_{k}, v\in X$ is continuous. Thus 
$\{g_{k}(\lambda\otimes {\cal X}_{\alpha})\}_{\alpha\in A}$ 
is convergent to $0$ in $F_{k}$. In particular, there exists 
$\alpha_{0}\in A$ such that 
$\{g_{k}(\lambda\otimes {\cal X}_{\alpha})\}_{\alpha\in A, 
\alpha>\alpha_{0}}$ is weakly bounded. Thus 
$$\left\{\sup_{\alpha'\in A, \alpha'>\alpha_{0}}\mu_{\alpha}(
g_{k}(\lambda\otimes {\cal X}_{\alpha'}))\right\}_{\alpha\in A}$$
is convergent to $0$. In particular,
$$\{\mu_{\alpha}(
g_{k}(\lambda\otimes {\cal X}_{\alpha}))\}_{\alpha\in A, \alpha>\alpha_{0}}
$$
or equivalently 
$$\{\mu_{\alpha}(
g_{k}(\lambda\otimes {\cal X}_{\alpha}))\}_{\alpha\in A}
$$
is convergent to $0$. Since
$$\langle \lambda, e_{k}({\cal X}_{\alpha}\otimes
\mu_{\alpha})
\rangle=\mu_{\alpha}(g_{k}(\lambda\otimes {\cal X}_{\alpha})),
$$
we see that $\{\langle \lambda, e_{k}({\cal X}_{\alpha}\otimes
\mu_{\alpha})
\rangle\}_{\alpha\in A}$ is convergent to $0$ for $\lambda\in \tilde{G}$.
By the definition of the topology on $\tilde{G}^{*}$,
$\{e_{k}({\cal X}_{\alpha}\otimes
\mu_{\alpha})\}_{\alpha\in A}$ is convergent to $0$ in 
$\tilde{G}^{*}$.
\epfv

From Proposition \ref{elcont}, we conclude:

\begin{corol}
The quotient space 
$$(X^{\otimes (k+1)}\otimes
F^{*}_{k})/(e_{k}|_{X^{\otimes (k+1)}\otimes F^{*}_{k}})
^{-1}(0)$$
is a  locally convex
space. \epf
\end{corol}

Since $G_{k}$ is linearly isomorphic to 
$$(X^{\otimes (k+1)}\otimes
F_{k}^{*})/(e_{k}|_{X^{\otimes (k+1)}\otimes F_{k}^{*}})^{-1}(0),$$
the  
locally convex
space structure on 
$$(X^{\otimes (k+1)}\otimes
F_{k}^{*})/(e_{k}|_{X^{\otimes (k+1)}\otimes F_{k}^{*}})^{-1}(0)$$
gives a locally convex
space structure on $G_{k}$. Let $H_{k}$ be the completion of
$G_{k}$. Then $H_{k}$ is a complete locally convex space.

\begin{propo}\label{ind-tp}
The space $H_{k}$ can be embedded canonically in $H_{k+1}$. 
The topology on $H_{k}$ is the same as the one induced from the
topology on $H_{k+1}$. 
\end{propo}
\pf
The first conclusion follows from Proposition \ref{inclusion}. 

To prove the second conclusion, we need only prove that the topology
on $G_{k}$ is the same as the one induced from the topology on
$G_{k+1}$. 
We have the following commutative diagram:
$$\begin{CD}
X^{\otimes (k+1)}\otimes F_{k}^{*}@>\psi_{k}>>X^{\otimes (k+2)}\otimes
F_{k+1}^{*}\\
@Ve_{k}VV @VVe_{k+1}V\\
G_{k}@>>>G_{k+1}
\end{CD}$$
where 
$$\psi_{k}: X^{\otimes (k+1)}\otimes F_{k}^{*}\to
X^{\otimes (k+2)}\otimes F_{k+1}^{*}$$
is defined by
$$\psi_{k}(u_{1}\otimes \cdots \otimes u_{k}\otimes v\otimes \mu)
={\bf 1}\otimes u_{1}\otimes \cdots \otimes u_{k}\otimes v\otimes
\gamma_{k}^{*}(\mu)$$
for $u_{1}, \dots, u_{k}, v\in X$ and $\mu\in F_{k}^{*}$. 
By Corollary \ref{k*:k+1*},  if we identify $F_{k}^{*}$ with 
$\gamma_{k}^{*}(F_{k}^{*})$,
the topology on $F_{k}^{*}$ is 
induced from the one on $F_{k+1}^{*}$. By the definition of the topologies
on $G_{k}$ and $G_{k+1}$ 
and the commutativity of the diagram above, we see that
the topology on $G_{k}$ is the same as the one induced from the
topology on $G_{k+1}$. 
\epfv

By Proposition \ref{ind-tp}, we have a sequence $\{H_{k}\}_{k\ge 0}$ 
of strictly
increasing complete locally convex
spaces.
Let 
$$H=\bigcup_{k\ge 0}H_{k}.$$
We equip $H$ with the inductive
limit topology. Then $H$ is a complete locally convex space. 
Let 
$$G=\bigcup_{k\ge 0}G_{k}.$$
 Then $G$ is a dense subspace of $H$. 
Note that $V\subset G
\subset \overline{V}$. Since  elements of $\overline{V}$ are
finite or infinite sums of elements of $V$, 
elements of $G$ are also
finite or infinite sums of elements of $V$. Thus infinite sums 
of elements of $V$ belonging to
$G$ must be convergent in the topology on $H$. So $G$ is in the
closure of $V$. Since $G$ is dense in $H$, we obtain:

\begin{theo}
The vector space $H$ equipped with the strict inductive limit topology
is a locally convex 
completion of $V$.\epf
\end{theo}

\renewcommand{\theequation}{\thesection.\arabic{equation}}
\renewcommand{\therema}{\thesection.\arabic{rema}}
\setcounter{equation}{0}
\setcounter{rema}{0}

\section{The locally convex completion and the vertex operator map}

In this section we construct continuous linear maps from 
the topological completion of $H\otimes H$
to $H$ associated with conformal equivalence classes of 
closed disks with two ordered open disks inside deleted and with the standard
parametrizations at the boundary components.

We consider the closed unit disk centered at $0$. We delete two
ordered open 
disks inside: the first
 centered at $z$ with radius $r_{1}$ and the second centered
at $0$ with radius $r_{2}$. See the figure below.

\begin{center}
\includegraphics[scale=.8]{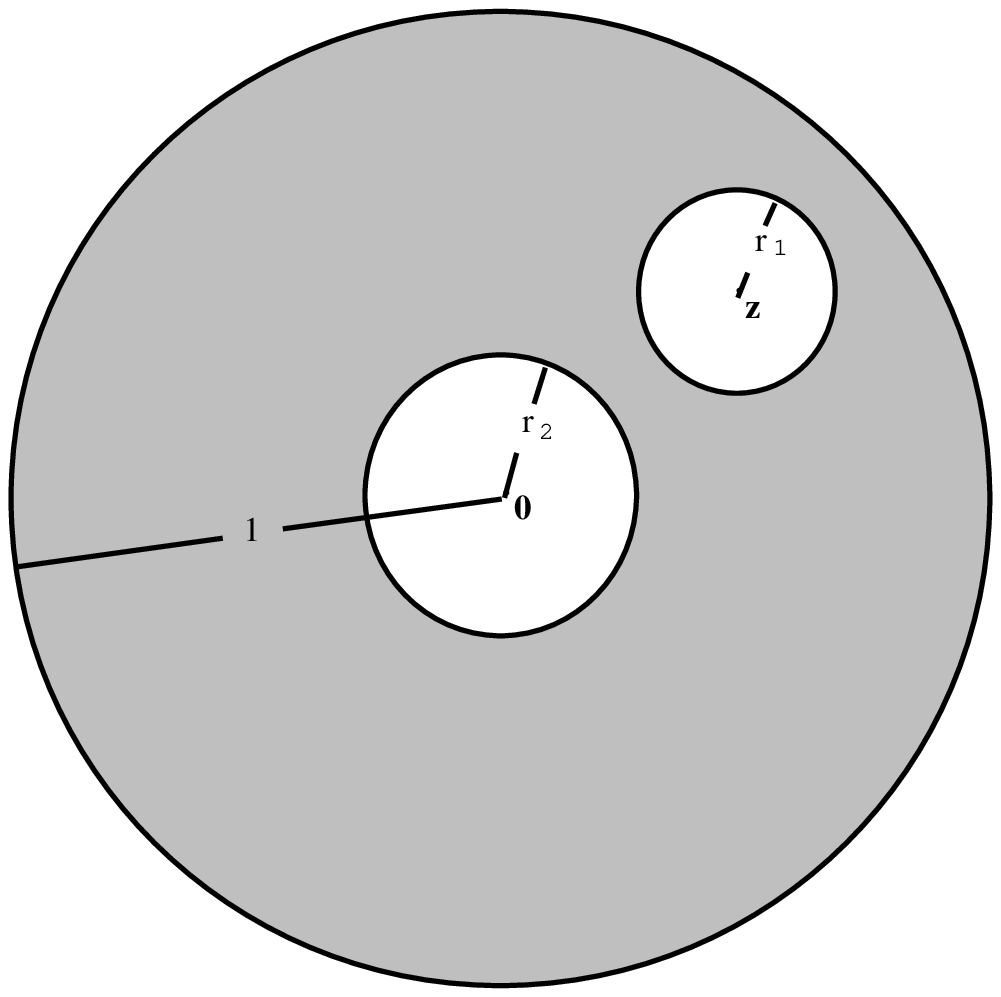}
\end{center}

The positive numbers $r_{1}, r_{2}$ and
the nonzero complex number $z$ must satisfy the conditions
$r_{2}+2r_{1}<1$ and $r_{2}<|z|<1$. The three boundary circles
are parametrized by the maps 
\begin{eqnarray*}
e^{i\theta}&\mapsto &e^{i\theta}\nno\\
e^{i\theta}&\mapsto &z+r_{1}e^{i\theta}\nno\\
e^{i\theta}&\mapsto &r_{2}e^{i\theta}.
\end{eqnarray*}
We denote the resulting closed 
disk with two disjoint open disks inside deleted and with ordered
boundary components parametrized as above by $D(z; r_{1},
r_{2})$. 
Note that any closed disks with two ordered open disks inside
deleted and with the standard  parametrizations is conformally
equivalent to $D(z; r_{1},
r_{2})$ for some $z\in {\Bbb C}^{\times}$, and some positive 
numbers $r_{1}, r_{2}$ satisfying $r_{2}+2r_{1}<1$ and $r_{2}<|z|<1$.

Let $H\widetilde{\otimes} H$ be 
the topological completion of the vector space tensor product
$H\otimes H$. We would like to construct a continuous linear map 
$$\overline{\nu}_{Y}([D(z;  r_{1},
r_{2})]): H\widetilde{\otimes} H\to H$$ 
associated with the conformal equivalence class $[D(z; r_{1},
r_{2})]$ of $D(z;  r_{1},
r_{2})$.  We know that $D(z;  r_{1},
r_{2})$ corresponds to the sphere $\hat{\Bbb C}=
{\Bbb C}\cup \{\infty\}$ with the negatively oriented 
puncture $\infty$, the ordered positively oriented puncture 
$z, 0$ and with the local coordinates $w\mapsto 1/w$, $w\mapsto 
(w-z)/r_{1}$ and 
$w/r_{2}$ vanishing at $\infty$, $z$ and $0$, respectively.
In the notation of \cite{H6}, 
the conformal equivalence class of this sphere with tubes of type $(1, 2)$
is 
$(z; {\bf 0}, (1/r_{1}, {\bf 0}), (1/r_{2}, {\bf 0}))\in K(2)$.
Associated with this class, we  have a linear map
$$\nu_{Y}((z; {\bf 0}, (1/r_{1}, {\bf 0}), (1/r_{2}, {\bf
0})))
=Y(r_{1}^{L(0)}\cdot, z)r_{2}^{L(0)}\cdot : 
V\otimes V
\to \overline{V}$$
(see \cite{H6}).
We now show that this linear map can be extended to
$H\widetilde{\otimes} H$,
that the
image of this extension 
is in $H$ and  that this extension is continuous. Then we 
define $\overline{\nu}_{Y}([D(z; r_{1}, r_{2})])$ to be this extension.

For any $\lambda\in \tilde{G}$ and $u\in V$, we define an element 
$u\diamond_{[D(z; r_{1}, r_{2})]}\lambda\in V^{*}$ by
\begin{eqnarray*}
(u\diamond_{[D(z; r_{1}, r_{2})]}\lambda)(v)
&
=&\sum_{n\in {\Bbb Z}}
\lambda(P_{n}(Y(r_{1}^{L(0)}u, z)r_{2}^{L(0)}v))\nno\\
&=&\sum_{n\in {\Bbb Z}}\lambda(P_{n}(Q(r_{1}^{L(0)}u, r_{2}^{L(0)}v; z)))
\end{eqnarray*}
for $v\in V$.

\begin{propo}\label{3-1}
The element $u\diamond_{[D(z; r_{1}, r_{2})]}\lambda$ is in $\tilde{G}$.
\end{propo}
\pf
For any $k\ge 0$, $u_{1}, \dots, u_{k},
v\in V$,
\begin{eqnarray}\label{3.1}
\lefteqn{\sum_{m\in {\Bbb Z}} 
(u\diamond_{[D(z; r_{1}, r_{2})]}\lambda)(P_{m}(Q(u_{1}, \dots, u_{k}, v;
z_{1}, \dots, z_{k})))}\nno\\ 
&&=\sum_{m\in {\Bbb Z}}\sum_{n\in
{\Bbb Z}} \lambda(P_{n}(Y(r_{1}^{L(0)}u, z)r_{2}^{L(0)}
P_{m}(Q(u_{1}, \dots, u_{k}, v;
z_{1}, \dots, z_{k})))\nno\\
&&=\sum_{m\in {\Bbb Z}}\sum_{n\in
{\Bbb Z}} \lambda(
P_{n}(Y(r_{1}^{L(0)}u, z)\cdot\nno\\
&&\hspace{5em}\cdot P_{m}(Q(r_{2}^{L(0)}u_{1}, \dots,
r_{2}^{L(0)}u_{k}, 
r_{2}^{L(0)}v;
r_{2}z_{1}, \dots, r_{2}z_{k}))).
\end{eqnarray}

Since the ${\Bbb Z}$-grading of $V$ is lower-truncated and 
since $r_{2}<|z|<1$, for any $(z_{1}, \dots, z_{k})\in M^{k}_{<1}$,
there exists a positive number $\delta_{2}>1$
such that
the Laurent series in $t_{2}$
\begin{eqnarray*}
\lefteqn{\sum_{m\in {\Bbb Z}}\lambda(
P_{n}(Y(r_{1}^{L(0)}u, z)P_{m}(Q(r_{2}^{L(0)}u_{1}, \dots,
r_{2}^{L(0)}u_{k}, 
r_{2}^{L(0)}v;
r_{2}z_{1}, \dots, r_{2}z_{k}))))t_{2}^{m}}\nno\\
&&=\sum_{m\in {\Bbb Z}}\lambda(
P_{n}(Y(r_{1}^{L(0)}u, z)\cdot\nno\\
&&\hspace{5em}\cdot P_{m}(t_{2}^{L(0)}Q(r_{2}^{L(0)}u_{1}, \dots,
r_{2}^{L(0)}u_{k}, 
r_{2}^{L(0)}v;
r_{2}z_{1}, \dots, r_{2}z_{k}))))\nno\\
&&=\sum_{m\in {\Bbb Z}}\lambda(
P_{n}(Y(r_{1}^{L(0)}u, z)P_{m}(Q((t_{2}r_{2})^{L(0)}u_{1}, \dots,
(t_{2}r_{2})^{L(0)}u_{k}, 
(t_{2}r_{2})^{L(0)}v;\nno\\
&&\hspace{5em} t_{2}r_{2}z_{1}, \dots, t_{2}r_{2}z_{k}))))
\end{eqnarray*}
in $t_{2}$ has only finitely many negative powers and 
 is absolutely convergent to 
$$\lambda(P_{n}(Q(r_{1}^{L(0)}u, (t_{2}r_{2})^{L(0)}u_{1}, \dots,
(t_{2}r_{2})^{L(0)}u_{k}, 
(t_{2}r_{2})^{L(0)}v; z, t_{2}r_{2}z_{1}, \dots, 
t_{2}r_{2}z_{k})))$$
when $0<|t_{2}|<\delta_{2}$.  
Since $\lambda\in \tilde{G}$
and $|z|<1$, for any $(z_{1}, \dots, z_{k})\in M^{k}_{<1}$, 
there exists a positive number $\delta_{1}>1$
such that the Laurent series in $t_{1}$
\begin{eqnarray*}
\lefteqn{\sum_{n\in
{\Bbb Z}} \lambda(
P_{n}(Q(r_{1}^{L(0)}u, r_{2}^{L(0)}u_{1}, \dots, r_{2}^{L(0)}u_{k}, 
r_{2}^{L(0)}v; z, r_{2}z_{1}, \dots, 
r_{2}z_{k})))t_{1}^{n}}\nno\\
&&=\sum_{n\in
{\Bbb Z}} \lambda(
P_{n}(t_{1}^{L(0)}
Q(r_{1}^{L(0)}u, r_{2}^{L(0)}u_{1}, \dots, r_{2}^{L(0)}u_{k}, 
r_{2}^{L(0)}v; z, r_{2}z_{1}, \dots, 
r_{2}z_{k})))\nno\\
&&=\sum_{n\in
{\Bbb Z}} \lambda(
P_{n}(Q((t_{1}r_{1})^{L(0)}u, 
(t_{1}r_{2})^{L(0)}u_{1}, \dots, (t_{1}r_{2})^{L(0)}u_{k}, 
(t_{1}r_{2})^{L(0)}v; \nno\\
&&\hspace{5em} t_{1}z, t_{1}r_{2}z_{1}, \dots, 
t_{1}r_{2}z_{k})))
\end{eqnarray*}
is absolutely convergent when $0<|t_{1}|<\delta_{1}$. 
Thus for $(z_{1}, \dots, z_{k})\in M^{k}_{<1}$, the iterated sum
\begin{eqnarray}\label{3.1.5}
&{\dps \sum_{n\in
{\Bbb Z}} \sum_{m\in {\Bbb Z}}
\lambda(
P_{n}(Y(r_{1}^{L(0)}u, z)P_{m}(Q(r_{2}^{L(0)}u_{1}, \dots,
r_{2}^{L(0)}u_{k}, 
r_{2}^{L(0)}v;}&\nno\\
&r_{2}z_{1}, \dots, r_{2}z_{k}))))t_{1}^{n}t_{2}^{m}&
\end{eqnarray}
is absolutely convergent when $0<|t_{1}|<\delta_{1}$ and 
$0<|t_{2}|<\delta_{2}$. 

The iterated sum of (\ref{3.1.5}) gives a function
of $t_{1}$ and $t_{2}$ in the region $0<|t_{1}|<\delta_{1}$,
$0<|t_{2}|<\delta_{2}$. It is clear that this function
is analytic in $t_{1}$ and $t_{2}$. Thus it has a Laurent 
expansion which must be 
\begin{eqnarray*}
&{\dps \sum_{n, m\in
{\Bbb Z}} 
\lambda(
P_{n}(Y(r_{1}^{L(0)}u, z)P_{m}(Q(r_{2}^{L(0)}u_{1}, \dots,
r_{2}^{L(0)}u_{k}, 
r_{2}^{L(0)}v;}\nno\\
& r_{2}z_{1}, \dots, r_{2}z_{k}))))t_{1}^{n}t_{2}^{m}.&
\end{eqnarray*}
Since this double sum is equal to the Laurent expansion, 
it is absolutely convergent and thus both iterated sums are 
absolutely convergent. In particular, when $t_{1}=t_{2}=1$, we see that
the right-hand side of (\ref{3.1}) and consequently the left-hand side
of (\ref{3.1}) are absolutely convergent 
when $(z_{1}, \dots, z_{k})\in M^{k}_{<1}$, 
proving that $u\diamond_{[D(z; r_{1}, r_{2})]}\lambda$ is in $\tilde{G}$.
\epfv

For any $l\ge 0$, 
we define a linear map 
$\alpha_{l}: \tilde{G} \otimes X^{l+1}\otimes F_{l}^{*}\to V^{*}$ by
\begin{eqnarray*}
\lefteqn{(\alpha_{l}(\lambda\otimes v_{1}\otimes \dots 
\otimes v_{l}\otimes v\otimes \nu))(u)}\nno\\
&&=\langle u\diamond_{[D(z; r_{1}, r_{2})]}\lambda,
e_{l}
(v_{1}\otimes \dots 
\otimes v_{l}\otimes v\otimes \nu)\rangle.
\end{eqnarray*}
for $\lambda\in \tilde{G}$,  $v_{1}, \dots, v_{l}, v\in X$,
$\nu\in F_{l}^{*}$ and $u\in V$.

\begin{propo}\label{3-2}
The image of  $\alpha_{l}$ is  in $\tilde{G}$.
\end{propo}
\pf
For any $k\ge 0$,
$\lambda\in \tilde{G}$, $u_{1}, \dots, u_{k}, u, \in V$, 
$v_{1}, \dots, v_{l},  v\in X$
and 
$\nu\in F_{l}^{*}$,
\begin{eqnarray}\label{3.2}
\lefteqn{\sum_{n\in {\Bbb Z}}
(\alpha_{l}(\lambda\otimes v_{1}\otimes \dots 
\otimes v_{l}\otimes v\otimes \nu))(P_{n}(Q(u_{1}, \dots, u_{k}, u; 
z_{1}, \dots, z_{k})))}\nno\\
&&=\sum_{n\in {\Bbb Z}}\langle (P_{n}(Q(u_{1}, \dots, u_{k}, u; 
z_{1}, \dots, z_{k}))\diamond_{[D(z; r_{1}, r_{2})]}\lambda),\nno\\
&&\hspace{6em}e_{l}
(v_{1}\otimes \cdots \otimes v_{l} \otimes
v\otimes \nu)\rangle\nno\\
&&=\sum_{n\in {\Bbb Z}}\nu(g_{l}((P_{n}(Q(u_{1}, \dots, u_{k}, u; 
z_{1}, \dots, z_{k}))\diamond_{[D(z; r_{1}, r_{2})]}\lambda)\nno\\
&&\hspace{6em}\otimes 
v_{1}\otimes
\cdots \otimes \otimes v_{l}\otimes v))\nno\\
&&=\sum_{n\in {\Bbb Z}}
\nu\Biggl(\sum_{m\in {\Bbb Z}}(P_{n}(Q(u_{1}, \dots, u_{k}, u; 
z_{1}, \dots, z_{k}))\diamond_{[D(z; r_{1}, r_{2})]}\lambda)\nno\\
&&\hspace{6em}(
P_{m}(Q(v_{1},
\dots, v_{k},
v; z_{k+1}, \dots, z_{k+l})))\Biggr).
\end{eqnarray}
By the definition of 
$$P_{n}(Q(u_{1}, \dots, u_{k}, u; 
z_{1}, \dots, z_{k}))\diamond_{[D(z; r_{1}, r_{2})]}\lambda,$$
we have 
\begin{eqnarray}\label{3.3}
\lefteqn{(P_{n}(Q(u_{1}, \dots, u_{k}, u; 
z_{1}, \dots, z_{k}))\diamond_{[D(z; r_{1}, r_{2})]}\lambda)}\nno\\
&&\hspace{6em}(P_{m}(Q(v_{1},
\dots, v_{k},
v; z_{k+1}, \dots, z_{k+l})))\nno\\
&&=\sum_{p\in {\Bbb Z}}\lambda(P_{p}(
Y(r_{1}^{L(0)}P_{n}(Q(u_{1}, \dots, u_{k}, u; 
z_{1}, \dots, z_{k})), z)\cdot \nno\\
&&\hspace{6em} \cdot r_{2}^{L(0)}P_{m}(Q(v_{1},
\dots, v_{l},
v; z_{k+1}, \dots, z_{k+l}))))\nno\\
&&=\sum_{p\in {\Bbb Z}}\lambda(P_{p}(
Y(P_{n}(Q(r_{1}^{L(0)}u_{1}, \dots, r_{1}^{L(0)}u_{k}, r_{1}^{L(0)}u; 
r_{1}z_{1}, \dots, r_{1}z_{k})), z)\cdot \nno\\
&&\hspace{6em}\cdot P_{m}(Q(r_{2}^{L(0)}v_{1},
\dots, r_{2}^{L(0)}v_{l},
r_{2}^{L(0)}v; r_{2}z_{k+1}, \dots, r_{2}z_{k+l})))).\nno\\
&&
\end{eqnarray}

Using the associativity of vertex operators, we know that the element
\begin{eqnarray*}
\lefteqn{\sum_{n\in {\Bbb Z}}\sum_{m\in {\Bbb Z}}
Y(P_{n}(Q(r_{1}^{L(0)}u_{1}, \dots, r_{1}^{L(0)}u_{k}, r_{1}^{L(0)}u; 
r_{1}z_{1}, \dots, r_{1}z_{k})), z) \cdot} \nno\\
&&\hspace{5em}\cdot P_{m}(Q(r_{2}^{L(0)}v_{1},
\dots, r_{2}^{L(0)}v_{l},
r_{2}^{L(0)}v; r_{2}z_{k+1}, \dots, r_{2}z_{k+l}))
\end{eqnarray*}
in $\overline{V}$ is  equal to 
\begin{eqnarray*}
\lefteqn{Q(r_{1}^{L(0)}u_{1}, \dots, r_{1}^{L(0)}u_{k}, r_{1}^{L(0)}u, 
r_{2}^{L(0)}v_{1},
\dots, r_{2}^{L(0)}v_{l},
r_{2}^{L(0)}v; }\nno\\
&&\hspace{5em}
r_{1}z_{1}+z, \dots, r_{1}z_{k}+z, z, r_{2}z_{k+1}, \dots,
r_{2}z_{k+l})
\end{eqnarray*}
when $(z_{1}, \dots, z_{k})\in
M^{k}_{<1}$ and $(z_{k+1}, \dots, z_{k+l})\in M^{l}_{<1}$. 
This fact and (\ref{3.3}) imply that
\begin{eqnarray*}
\lefteqn{\sum_{n\in {\Bbb Z}}\sum_{m\in {\Bbb Z}}(P_{n}(Q(u_{1}, 
\dots, u_{k}, u; 
z_{1}, \dots, z_{k}))\diamond_{[D(z; r_{1}, r_{2})]}\lambda)}\nno\\
&&\hspace{5em}(
P_{m}(Q(v_{1},
\dots, v_{k}, v; z_{k+1}, \dots, z_{k+l})))
\end{eqnarray*}
is convergent absolutely to 
\begin{eqnarray*}
\lefteqn{\lambda(P_{p}(Q(r_{1}^{L(0)}u_{1}, 
\dots, r_{1}^{L(0)}u_{k}, r_{1}^{L(0)}u, 
r_{2}^{L(0)}v_{1},
\dots, r_{2}^{L(0)}v_{l},
r_{2}^{L(0)}v; }\nno\\
&&\hspace{5em}
r_{1}z_{1}+z, \dots, r_{1}z_{k}+z, z, r_{2}z_{k+1}, \dots,
r_{2}z_{k+l})))
\end{eqnarray*}
when $(z_{1}, \dots, z_{k})\in
M^{k}_{<1}$ and $(z_{k+1}, \dots, z_{k+l})\in M^{l}_{<1}$. 
Since  $\lambda\in \tilde{G}$,
\begin{eqnarray*}
\lefteqn{\sum_{p\in {\Bbb Z}}\lambda(P_{p}(Q(r_{1}^{L(0)}u_{1}, 
\dots, r_{1}^{L(0)}u_{k}, r_{1}^{L(0)}u, 
r_{2}^{L(0)}v_{1},
\dots, r_{2}^{L(0)}v_{l},
r_{2}^{L(0)}v; }\nno\\
&&\hspace{5em}
r_{1}z_{1}+z, \dots, r_{1}z_{k}+z, z, r_{2}z_{k+1}, \dots,
r_{2}z_{k+l})))
\end{eqnarray*}
is absolutely convergent when $(z_{1}, \dots, z_{k})\in
M^{k}_{<1}$ and $(z_{k+1}, \dots, z_{k+l})\in M^{l}_{<1}$.
The arguments above prove
that the iterated sum
\begin{eqnarray*}
\lefteqn{\sum_{p\in {\Bbb Z}}\sum_{n\in {\Bbb Z}}
\sum_{m\in {\Bbb Z}}\lambda(P_{p}(
Y(P_{n}(Q(r_{1}^{L(0)}u_{1}, \dots, r_{1}^{L(0)}u_{k}, r_{1}^{L(0)}u; 
r_{1}z_{1}, \dots, r_{1}z_{k})), z)\cdot }\nno\\
&&\hspace{5em} \cdot P_{m}(Q(r_{2}^{L(0)}v_{1},
\dots, r_{2}^{L(0)}v_{l},
r_{2}^{L(0)}v;  r_{2}z_{k+1}, \dots, r_{2}z_{k+l}))))
\end{eqnarray*}
is absolutely convergent 
when $(z_{1}, \dots, z_{k})\in
M^{k}_{<1}$ and $(z_{k+1}, \dots, z_{k+l})\in M^{l}_{<1}$.  
The same method as in the proof of Proposition \ref{3-1}
shows that in fact the corresponding triple sum is absolutely convergent
and thus all the iterated sums are absolutely convergent and are all
equal when $(z_{1}, \dots, z_{k})\in
M^{k}_{<1}$ and $(z_{k+1}, \dots, z_{k+l})\in M^{l}_{<1}$. 
So
\begin{eqnarray*}
\lefteqn{\sum_{n\in {\Bbb Z}}
\sum_{m\in {\Bbb Z}}\sum_{p\in {\Bbb Z}}\lambda(P_{p}(
Y(P_{n}(Q(r_{1}^{L(0)}u_{1}, \dots, r_{1}^{L(0)}u_{k}, r_{1}^{L(0)}u; 
r_{1}z_{1}, \dots, r_{1}z_{k})), z)\cdot }\nno\\
&&\hspace{5em}\cdot P_{m}(Q(r_{2}^{L(0)}v_{1},
\dots, r_{2}^{L(0)}v_{l},
r_{2}^{L(0)}v; r_{2}z_{k+1}, \dots, r_{2}z_{k+l}))))\nno\\
&&=\sum_{n\in {\Bbb Z}}
\sum_{m\in {\Bbb Z}}(P_{n}(Q(u_{1}, \dots, u_{k}, u; 
z_{1}, \dots, z_{k}))\diamond_{[D(z; r_{1}, r_{2})]}\lambda)\nno\\
&&\hspace{5em}(
P_{m}(Q(v_{1},
\dots, v_{k},
v; z_{k+1}, \dots, z_{k+l})))
\end{eqnarray*}
is absolutely convergent when $(z_{1}, \dots, z_{k})\in
M^{k}_{<1}$ and $(z_{k+1}, \dots, z_{k+l})\in M^{l}_{<1}$. 
By (\ref{3.2}), we see that the left-hand side of 
(\ref{3.2}) is also absolutely convergent when $(z_{1}, \dots, z_{k})\in
M^{k}_{<1}$ and $(z_{k+1}, \dots, z_{k+l})\in M^{l}_{<1}$, proving that 
$$\alpha_{l}(\lambda \otimes v_{1}\otimes \dots 
\otimes v_{l}\otimes v\otimes \nu)$$
is indeed in $\tilde{G}$.
\epfv

By Proposition \ref{3-2}, 
$$\sum_{n\in {\Bbb Z}}\alpha_{l}(\lambda \otimes v_{1}\otimes \dots 
\otimes v_{l}\otimes v\otimes \nu)(P_{n}(Q(u_{1}, \dots, u_{k}, u; 
z_{1}, \dots, z_{k})))$$
is absolutely convergent and the sum is equal to 
$$g_{k}(\alpha_{l}(\lambda\otimes v_{1}\otimes \dots 
\otimes v_{l}\otimes v\otimes \nu)\otimes u_{1}\otimes \cdots\otimes 
u_{k}\otimes u)\in F_{k}.$$
We define an linear map 
$$\beta_{k, l}: F_{k}^{*}\otimes F_{l}^{*}
\to F_{k+l+1}^{*}$$ 
by
\begin{eqnarray*}
\lefteqn{(\beta_{k, l}(\mu, \nu))(g_{k+l+1}(\lambda
\otimes 
u_{1}\otimes \dots \otimes u_{k+1}\otimes u\otimes 
v_{1}\otimes \cdots\otimes v_{l}\otimes v))}\nno\\
&&=\mu
(g_{k}(\alpha_{l}(\lambda \otimes v_{1}\otimes \dots 
\otimes v_{l}\otimes v\otimes \nu)\otimes u_{1}\otimes \cdots\otimes 
u_{k}\otimes u))
\end{eqnarray*}
for $\lambda\in \tilde{G}$, $u_{1}, \dots,  u_{k}, u, v_{1}, \dots,
v_{l}, v\in V$, $\mu\in F_{k}^{*}$ and $\nu\in F_{l}^{*}$. In fact
this formula only gives a linear map from $F_{k}^{*}\otimes
F_{l}^{*}$ to the algebraic dual of $F_{k+l+1}$.
We have:

\begin{propo}\label{3-4}
The image of the map $\beta_{k, l}$  is indeed in $F_{k+l+1}^{*}$ and the map
$\beta_{k, l}$ is continuous.
\end{propo}
\pf
To avoid notational confusions in the proof, we use 
$F_{k}^{\zeta}$, $k\ge 0$, to  denote the space $F_{k}$ 
whose elements  are viewed as
functions of $(\zeta_{1}, \dots, \zeta_{k})\in M^{k}_{<1}$. Similarly
we have the notation $F_{l}^{\eta}$, $l\ge 0$.
Let $\mu\in F_{k}^{*}$ and $\nu\in F_{l}^{*}$. 
For any $\lambda\in \tilde{G}$, $u_{1}, \dots,  u_{k}, u$, $v_{1}, \dots,
v_{l}, v\in V$, we have an element 
\begin{equation}\label{3.3.9}
g(z_{1}, \dots, z_{k+l+1})=g_{k+l+1}(\lambda
\otimes 
u_{1}\otimes \dots \otimes u_{k+1}\otimes u\otimes 
v_{1}\otimes \cdots\otimes v_{l}\otimes v)
\end{equation}
of $F_{k+l+1}$.
By definition 
\begin{eqnarray}\label{3.4}
\lefteqn{|(\beta_{k, l}(\mu, \nu))(g(z_{1}, \dots, z_{k+l+1}))|}\nno\\
&&=|(\beta_{k, l}(\mu, \nu))(g_{k+l+1}(\lambda
\otimes 
u_{1}\otimes \dots \otimes u_{k+1}\otimes u\otimes 
v_{1}\otimes \cdots\otimes v_{l}\otimes v))|\nno\\
&&=|\mu
(g_{k}(\alpha_{l}(\lambda \otimes v_{1}\otimes \dots 
\otimes v_{l}\otimes v\otimes \nu)\otimes u_{1}\otimes \cdots\otimes 
u_{k}\otimes u))|.
\end{eqnarray}
We now view
$$g_{k}(\alpha_{l}(\lambda \otimes v_{1}\otimes \dots 
\otimes v_{l}\otimes v\otimes \nu)\otimes u_{1}\otimes \cdots\otimes 
u_{k}\otimes u)$$ 
as an element of $F_{k}^{\zeta}$. 
Then 
\begin{eqnarray}\label{3.5}
\lefteqn{g_{k}(\alpha_{l}(\lambda \otimes v_{1}\otimes \dots 
\otimes v_{l}\otimes v\otimes \nu)\otimes u_{1}\otimes \cdots\otimes 
u_{k}\otimes u)}\nno\\
&&= \sum_{p\in {\Bbb Z}}(\alpha_{l}(\lambda \otimes v_{1}\otimes \dots 
\otimes v_{l}\otimes v\otimes \nu))(
P_{p}(Q(u_{1}, \dots, u_{k}, u; \zeta_{1}, \dots,
\zeta_{k})))\nno\\
&&=\sum_{p\in {\Bbb Z}}\nu(g_{l}((P_{p}(Q(u_{1}, \dots, u_{k}, 
u; \zeta_{1}, \dots,
\zeta_{k}))\diamond_{[D(z; r_{1}, r_{2})]}\lambda)\nno\\
&&\hspace{6em}
\otimes v_{1}\otimes \cdots \otimes 
v_{l}\otimes v)).
\end{eqnarray}

For fixed $\zeta_{1}, \dots, \zeta_{k}$, we view
$$g_{l}((P_{p}(Q(u_{1}, \dots, u_{k}, 
u; \zeta_{1}, \dots,
\zeta_{k_{1}}))\diamond_{[D(z; r_{1}, r_{2})]}\lambda)
\otimes v_{1}\otimes \cdots \otimes 
v_{l}\otimes v)$$
as an element of $F_{l}^{\eta}$. 
Then 
by definition we have
\begin{eqnarray}\label{3.6}
\lefteqn{\sum_{p\in {\Bbb Z}}g_{l}((P_{p}(Q(u_{1}, \dots, u_{k}, 
u; \zeta_{1}, \dots,
\zeta_{k_{1}}))\diamond_{[D(z; r_{1}, r_{2})]}\lambda)
\otimes v_{1}\otimes \cdots \otimes 
v_{l}\otimes v)}\nno\\
&&=\sum_{p\in {\Bbb Z}}\sum_{q\in {\Bbb Z}}
(P_{p}(Q(u_{1}, \dots, u_{k}, 
u; \zeta_{1}, \dots,
\zeta_{k_{1}}))\diamond_{[D(z; r_{1}, r_{2})]}\lambda)\nno\\
&&\hspace{6em}
(P_{q}(Q(v_{1}, \dots,
v_{l},  v; \eta_{1}, \dots, \eta_{l})))\nno\\
&&=\sum_{p\in {\Bbb Z}}\sum_{q\in {\Bbb Z}}\sum_{j\in {\Bbb Z}}
\lambda(P_{j}(Y(r_{1}^{L(0)}P_{p}(Q(u_{1}, \dots, u_{k}, 
u; \zeta_{1}, \dots,
\zeta_{k_{1}})), z) \cdot\nno\\
&&\hspace{6em}\cdot r_{2}^{L(0)}P_{q}(Q(v_{1}, \dots,
v_{l},  v; \eta_{1}, \dots, \eta_{l}))))\nno\\
&&=\sum_{p\in {\Bbb Z}}\sum_{q\in {\Bbb Z}}\sum_{j\in {\Bbb Z}}
\lambda(P_{j}(Y(P_{p}(Q(r_{1}^{L(0)}u_{1}, \dots, r_{1}^{L(0)}u_{k}, 
r_{1}^{L(0)}u; r_{1}\zeta_{1}, \dots,
r_{1}\zeta_{k_{1}})), z) \cdot \nno\\
&&\hspace{6em} \cdot P_{q}(Q(r_{2}^{L(0)}v_{1}, \dots,
r_{2}^{L(0)}v_{l},  r_{2}^{L(0)}v; r_{2}\eta_{1}, \dots, 
r_{2}\eta_{l})))).
\end{eqnarray}

We now calculate the right-hand side of (\ref{3.6}).
For any nonzero complex numbers $t_{0}, t_{1}, t_{2}$,
\begin{eqnarray}\label{3.6.1}
\lefteqn{\sum_{p, q, j\in {\Bbb Z}}
\lambda(P_{j}(Y(P_{p}(Q(r_{1}^{L(0)}u_{1}, \dots, r_{1}^{L(0)}u_{k}, 
r_{1}^{L(0)}u; r_{1}\zeta_{1}, \dots,
r_{1}\zeta_{k_{1}})), z) \cdot} \nno\\
&&\hspace{6em} \cdot P_{q}(Q(r_{2}^{L(0)}v_{1}, \dots,
r_{2}^{L(0)}v_{l},  r_{2}^{L(0)}v; r_{2}\eta_{1}, \dots, 
r_{2}\eta_{l}))))t_{0}^{j}t_{1}^{p}t_{2}^{q}\nno\\
&&=\sum_{p, q, j\in {\Bbb Z}}
\lambda(P_{j}(t_{0}^{L(0)}Y(P_{p}(t_{1}^{L(0)}
Q(r_{1}^{L(0)}u_{1}, \dots, r_{1}^{L(0)}u_{k}, 
r_{1}^{L(0)}u; \nno\\
&&\hspace{8em} r_{1}\zeta_{1}, \dots,
r_{1}\zeta_{k_{1}})), z) \cdot \nno\\
&&\hspace{6em} \cdot P_{q}(t_{2}^{L(0)}Q(r_{2}^{L(0)}v_{1}, \dots,
r_{2}^{L(0)}v_{l},  r_{2}^{L(0)}v; r_{2}\eta_{1}, \dots, 
r_{2}\eta_{l}))))\nno\\
&&=\sum_{p, q, j\in {\Bbb Z}}
\lambda(P_{j}(Y(P_{p}(
Q((t_{0}t_{1}r_{1})^{L(0)}u_{1}, \dots, 
(t_{0}t_{1}r_{1})^{L(0)}u_{k}, 
(t_{0}t_{1}r_{1})^{L(0)}u; \nno\\
&&\hspace{8em} t_{0}t_{1}r_{1}\zeta_{1}, \dots,
t_{0}t_{1}r_{1}\zeta_{k_{1}})), t_{0}z) \cdot \nno\\
&&\hspace{6em} \cdot P_{q}(Q((t_{0}t_{2}r_{1})^{L(0)}v_{1}, \dots,
(t_{0}t_{2}r_{1})^{L(0)}v_{l},  
(t_{0}t_{2}r_{1})^{L(0)}v; \nno\\
&&\hspace{8em}
t_{0}t_{2}r_{2}\eta_{1}, \dots, 
t_{0}t_{2}r_{2}\eta_{l})))).
\end{eqnarray}
By the associativity of vertex operators and the definition of the 
map $g_{k+l+1}$, there exists real numbers $\delta_{0},
\delta_{1}, \delta_{2}>1$ such that 
the right-hand side of (\ref{3.6.1}) is convergent absolutely 
to 
\begin{eqnarray}\label{3.6.2}
\lefteqn{g_{k+l+1}(\lambda \otimes (t_{0}t_{1}r_{1})^{L(0)}u_{1}\otimes 
\cdots \otimes (t_{0}t_{1}r_{1})^{L(0)}u_{k}\otimes (t_{0}t_{1}r_{1})^{L(0)}u
\otimes  }\nno\\
&&\hspace{2em} (t_{0}t_{2}r_{1})^{L(0)}v_{1}\otimes  \cdots \otimes 
(t_{0}t_{2}r_{1})^{L(0)}v_{l}\otimes \nno\\
&&\hspace{2em} (t_{0}t_{2}r_{1})^{L(0)}v)
\mbar_{z_{i}=t_{0}z+t_{0}t_{1}r_{1}\zeta_{i}, i=1,\dots, k,
z_{k+1}=t_{0}z,
z_{k+1+j}=t_{0}t_{2}r_{2}\eta_{j}, j=1, \dots, l}
\nno\\
&&
\end{eqnarray}
when $0<t_{0}<\delta_{0}$, $0<t_{1}<\delta_{1}$ and $0<t_{1}<\delta_{1}$.
In particular, when $t_{0}=t_{1}=t_{2}=1$, we see that (\ref{3.6.2}) 
is equal to
and thus the right-hand side of (\ref{3.6}) is convergent absolutely
to 
\begin{eqnarray}\label{3.6.3}
\lefteqn{g_{k+l+1}(\lambda \otimes r_{1}^{L(0)}u_{1}\otimes 
\cdots \otimes r_{1}^{L(0)}u_{k}\otimes r_{1}^{L(0)}u
\otimes r_{2}^{L(0)}v_{1}\otimes  \cdots}\nno\\
&&\hspace{2em} \otimes 
r_{2}^{L(0)}v_{l} \otimes r_{2}^{L(0)}v)
\mbar_{z_{i}=z+r_{1}\zeta_{i}, i=1,\dots, k,
 z_{k+1}=z,
z_{k+1+j}=r_{2}\eta_{j}, j=1, \dots, l}.
\nno\\
&&
\end{eqnarray}
Hence for fixed 
$(\zeta_{1}, \dots, \zeta_{k})\in M^{k}_{<1}$, 
 (\ref{3.6.3}) as a 
function of $\eta_{1}, \dots, \eta_{l}$ is an element of
$F_{l}^{\eta}$.
We denote this element 
by $f_{g;\zeta_{1}, \dots, \zeta_{k}}(\eta_{1}, \dots, \eta_{l})$.

We now view $\mu$ and
$\nu$ as elements of $(F_{k}^{\zeta})^{*}$ and $(F_{l}^{\eta})^{*}$,
respectively.
Since $\nu$ is continuous, by the calculations above,
we have 
\begin{eqnarray}\label{3.7}
\lefteqn{\nu(f_{g;\zeta_{1}, \dots, \zeta_{k}}(\eta_{1}, \dots, \eta_{l})
)}\nno\\
&&=\nu\Biggl(\sum_{p\in {\Bbb Z}}g_{l}((P_{p}(Q(u_{1}, \dots, u_{k}, 
u; \zeta_{1}, \dots,
\zeta_{k_{1}}))\diamond_{[D(z; r_{1}, r_{2})]}\lambda)\nno\\
&&\hspace{6em}
\otimes v_{1}\otimes \cdots \otimes 
v_{l}\otimes v)\Biggr)\nno\\
&&=\sum_{p\in {\Bbb Z}}\nu(g_{l}((P_{p}(Q(u_{1}, \dots, u_{k}, 
u; \zeta_{1}, \dots,
\zeta_{k_{1}}))\diamond_{[D(z; r_{1}, r_{2})]}\lambda)\nno\\
&&\hspace{6em}
\otimes v_{1}\otimes \cdots \otimes 
v_{l}\otimes v)).
\end{eqnarray}
By the calculations from (\ref{3.4}) to (\ref{3.7}), we see that
$$\nu(f_{g;\zeta_{1}, \dots, \zeta_{k}}(\eta_{1}, \dots, \eta_{l})
)\in F_{k}^{\zeta}$$
and
we obtain
\begin{eqnarray}\label{3.8}
\lefteqn{|(\beta_{k, l}(\mu, \nu))(g(z_{1}, \dots, z_{k+l+1}))|}\nno\\
&&=|\mu(\nu(f_{g;\zeta_{1}, \dots, \zeta_{k}}(\eta_{1}, \dots,
\eta_{l})))|.
\end{eqnarray}
For  $g(z_{1}, \dots, z_{k+l+1})\in 
F_{k+l+1}$ not of the form (\ref{3.3.9}), we use
linearity to define the functions $f_{g;\zeta_{1}, 
\dots, \zeta_{k}}(\eta_{1}, \dots, \eta_{l})$. Then
(\ref{3.8}) holds for any  $g(z_{1}, \dots, z_{k+l+1})\in 
F_{k+l+1}$. 

From the definition of $f_{g;\zeta_{1}, \dots, \zeta_{k}}(\eta_{1}, \dots,
\eta_{l})$, we see that for any
fixed $(\zeta_{1}, \dots, \zeta_{k})\in M^{k}_{<1}$, the linear
map from $F_{k+l+1}$ to $F_{l}^{\eta}$ given by
$$g\mapsto f_{g;\zeta_{1}, \dots, \zeta_{k}}(\eta_{1}, \dots, \eta_{l})$$
for $g\in F_{k+l+1}$ is continuous, and for any fixed 
$\nu\in (F_{l}^{\eta})^{*}$, 
the linear map from $F_{k+l+1}$ to $F_{k}^{\zeta}$ given by 
$$g\mapsto \nu(f_{g;\zeta_{1}, \dots, \zeta_{k}}(\eta_{1}, \dots, \eta_{l}))$$
is also continuous. Thus by (\ref{3.8}), $\beta_{k, l}(\mu, \nu)$ is
continuous. So the elements of image of $\beta_{k, l}$ are in 
$F_{k+l+1}^{*}$.

From the definition of 
$f_{g;\zeta_{1}, \dots, \zeta_{k}}(\eta_{1}, \dots, \eta_{l})$, we see that 
given any weakly bounded subset $B$ of $F_{k+l+1}$, 
the subset 
$$B'=\{f_{g;\zeta_{1}, \dots, \zeta_{k}}(\eta_{1}, \dots, \eta_{l})\;|\; 
g\in B\}\subset F_{l}^{\eta}$$
for fixed $(\zeta_{1}, \dots, \zeta_{k})\in M^{k}_{<1}$ is 
weakly bounded, and thus
for any net $\{\nu_{\alpha}\}_{\alpha\in A}$ convergent 
to $0$ in 
$(F_{l}^{\eta})^{*}$, there exists $\alpha_{0}\in A$ such that
the subset 
$$B''(\{\nu_{\alpha}\}_{\alpha\in A})=\{\nu_{\alpha}(f_{g;\zeta_{1}, 
\dots, \zeta_{k}}(\eta_{1}, \dots, \eta_{l}))\;|\; 
g\in B, \alpha>\alpha_{0}\}\subset F_{k}^{\zeta}$$
is also weakly bounded. 

Now let $\{(\mu_{\alpha}, \nu_{\alpha})\}_{\alpha\in A}$ be a 
net convergent to $0$ in $(F_{k}^{\zeta})^{*}\times (F_{l}^{\eta})^{*}$. 
Then the nets $\{\mu_{\alpha}\}_{\alpha\in A}$ and 
$\{\nu_{\alpha}\}_{\alpha\in A}$ 
are convergent to $0$ in $(F_{k}^{\zeta})$ and $(F_{l}^{\eta})$,
respectively. Thus by (\ref{3.8}) and the discussion above,
\begin{eqnarray}\label{3.9}
&\sup_{g(z_{1}, \dots, z_{k+l+1})\in B}
|(\beta_{k, l}(\mu_{\alpha}, \nu_{\alpha}))
(g(z_{1}, \dots, z_{k+l+1}))|&\nno\\
&\le \sup_{h(\zeta_{1}, \dots, \zeta_{k})\in 
B''(\{\nu_{\alpha}\}_{\alpha\in A})}
|\mu_{\alpha}(h(\zeta_{1}, \dots, \zeta_{k}))|&
\end{eqnarray}
when $\alpha>\alpha_{0}$. By the definition of the 
topology on $F_{k}^{*}$, the net on the right-hand side 
of (\ref{3.9}) is convergent to $0$. Thus the net on the left-hand
side of (\ref{3.9}) is convergent to $0$, proving the 
continuity of the map  $\beta_{k_{1}, k_{2}}$.
\epfv

Let 
$$h_{1}=e_{k}(u_{1}\otimes \cdots \otimes u_{k}\otimes
u\otimes \mu)\in G_{k}$$
and 
$$h_{2}=e_{l}(v_{1}\otimes \cdots \otimes v_{l}\otimes
v\otimes \nu)\in G_{l}$$
where $u_{1}, \dots, u_{k}, u, v_{1}, \dots, 
v_{l}, v\in X$, 
$\mu\in F_{k}^{*}$ and $\mu\in 
F_{l}^{*}$.
We define 
\begin{eqnarray*}
\lefteqn{(\overline{\nu}_{Y}([D(z; r_{0}, r_{1},
r_{2})]))(h_{1}\otimes h_{2})}\nno\\
&&=e_{k+l+1}(u_{1}\otimes \cdots \otimes 
u_{k}\otimes u\otimes v_{1}\otimes \cdots \otimes 
v_{l}\otimes  v\otimes \beta_{k, l}(\mu, \nu)).
\end{eqnarray*}
Note that any element of $G_{k}$ or $G_{l}$ is a linear
combination of elements of the form
$h_{1}$ or $h_{2}$, respectively,  given above, and that $k$ and $l$ are
arbitrary. Thus we obtain a linear map 
$$\overline{\nu}_{Y}([D(z; r_{0}, r_{1},
r_{2})])\mbar_{G\otimes G}: G\otimes G\to G.$$

\begin{propo}\label{contong}
The map $\overline{\nu}_{Y}([D(z; r_{1},
r_{2})])\mbar_{G\otimes G}$ is continuous.
\end{propo}
\pf
From the definition of $\overline{\nu}_{Y}([D(z; r_{1},
r_{2})])\mbar_{G\otimes G}$, we see that for any $k, l\ge 0$,
$$\overline{\nu}_{Y}([D(z; r_{1},
r_{2})]))(G_{k}\otimes G_{l})$$
is in $G_{k+l+1}$. 
To prove that $\overline{\nu}_{Y}([D(z; r_{1},
r_{2})])\mbar_{G\otimes G}$ is continuous, 
we need only prove that for any $k, l\ge 0$,
it is continuous as a map from $G_{k}\otimes G_{l}$ to $G_{k+l+1}$. 
Since the topology on $G_{k}$ and $G_{l}$ 
are defined to be  the quotient topologies on 
$$(X^{\otimes (k+1)}\otimes
F^{*}_{k})/(e_{k}|_{X^{\otimes (k+1)}\otimes
F^{*}_{k}})^{-1}(0)$$ and 
$$(X^{\otimes (l+1)}\otimes
F^{*}_{l})/(e_{l}|_{X^{\otimes (l+1)}\otimes
F^{*}_{l}})^{-1}(0),$$
respectively, we need only prove that the 
map 
\begin{eqnarray*}
\lefteqn{(\overline{\nu}_{Y}([D(z; r_{1},
r_{2})]))\circ(e_{k}|_{X^{\otimes (k+1)}\otimes
F^{*}_{k}}\otimes e_{l}|_{X^{\otimes (l+1)}\otimes
F^{*}_{l}}):}\nno\\
&&\hspace{3em} (X^{\otimes (k+1)}\otimes
F^{*}_{k})\otimes (X^{\otimes (l+1)}\otimes
F^{*}_{l})\to G_{k+l+1}
\end{eqnarray*}
is continuous. On the other hand, 
by definition, 
\begin{eqnarray*}
\lefteqn{(\overline{\nu}_{Y}([D(z; r_{1},
r_{2})]))\circ(e_{k}|_{X^{\otimes (k+1)}\otimes
F^{*}_{k}}\otimes e_{l}|_{X^{\otimes (l+1)}\otimes
F^{*}_{l}})}\nno\\ 
&&=e_{k+l+1}|_{X^{\otimes (k+l+2)}\otimes
F^{*}_{k+l+1}}\circ (I_{X^{\otimes (k+l+1)}}\otimes \beta_{k, l})\circ
\sigma_{23},
\end{eqnarray*}
where $I_{X^{\otimes (k+l+1)}}$ is the identity map on 
$X^{\otimes (k+l+1)}$ and 
$$\sigma_{23}: (X^{\otimes (k+1)}\otimes
F^{*}_{k})\otimes (X^{\otimes (l+1)}\otimes
F^{*}_{l})\to X^{\otimes (k+l+2)}\otimes F^{*}_{k}\otimes F^{*}_{l}$$ 
is the map permuting the second and the third tensor factor, that is,
$$\sigma_{23}({\cal X}_{1}\otimes \mu\otimes {\cal X}_{2}\otimes \nu)
={\cal X}_{1}\otimes {\cal X}_{2}\otimes  \mu\otimes \nu$$
for ${\cal X}_{1}\in X^{\otimes (k+1)}$, $\mu\in F^{*}_{k}$, 
${\cal X}_{2}\in X^{\otimes (l+1)}$ and $\nu\in F^{*}_{l}$.
Since the topology on $G_{k+l+1}$ is defined to be the quotient 
topology on 
$$(X^{\otimes (k+l+2)}\otimes
F^{*}_{k+l+1})/(e_{k+l+1}|_{X^{\otimes (k+l+2)}\otimes
F^{*}_{k+l+1}})^{-1}(0),$$
we need only prove that the map 
$(I_{X^{\otimes (k+l+1)}}\otimes \beta_{k, l})\circ
\sigma_{23}$ is continuous. But $I_{X^{\otimes (k+l+1)}}$ and 
$\sigma_{23}$ are obviously continuous and $\beta_{k, l}$ is continuous
by Proposition \ref{3-4}. So $(I_{X^{\otimes (k+l+1)}}\otimes 
\beta_{k, l})\circ
\sigma_{23}$ is indeed continuous. 
\epfv

Since $G$ is dense in $H$, we can extend $(\overline{\nu}_{Y}([D(z; r_{1},
r_{2})]))\mbar_{G\otimes G}$ to a linear map 
$\overline{\nu}_{Y}([D(z; r_{1},
r_{2})])$ from $H\widetilde{\otimes} H$ to $H$. 

\begin{theo}
The map $\overline{\nu}_{Y}([D(z; r_{1},
r_{2})])$ is a continuous extension of 
$$\nu_{Y}((z; {\bf 0}, (1/r_{1}, {\bf 0}), (1/r_{2}, {\bf
0})))$$ 
to $H\widetilde{\otimes}H.$
That is, $\overline{\nu}_{Y}([D(z; r_{1},
r_{2})])$ is continuous and
$$\overline{\nu}_{Y}([D(z; r_{1}, r_{2})])\mbar_{V\otimes
V}=\nu_{Y}((z; {\bf 0}, (1/r_{1}, {\bf 0}), (1/r_{2}, {\bf 0}))). $$
\end{theo}
\pf
The continuity of $\overline{\nu}_{Y}([D(z; r_{1},
r_{2})])$ follows from the definition  and Proposition \ref{contong}.

For any  $m_{1}, \dots, m_{k}\in {\Bbb Z}$,
let $\mu_{m_{1}, \dots, m_{k}}$ be an element of $F^{*}_{k}$ defined by
\begin{eqnarray*}
\lefteqn{\mu_{m_{1}, \dots, m_{k}}
(g_{k}(\lambda \otimes u_{1}\otimes \cdots \otimes u_{k}\otimes u))}\nno\\
&&=\frac{1}{2\pi \sqrt{-1}}\oint_{|z_{1}|=\epsilon_{1}}\cdots
\frac{1}{2\pi 
\sqrt{-1}}
\oint_{|z_{k}|=\epsilon_{k}}z_{1}^{m_{1}}\cdots 
z_{k}^{m_{k}}\cdot\nno\\
&&\hspace{6em} \cdot g_{k}(
\lambda \otimes u_{1}\otimes \cdots \otimes u_{k}\otimes u)dz_{k}\cdots 
dz_{k}
\end{eqnarray*}
for $\lambda\in \tilde{G}$, $u_{1}, \dots, u_{k}, u\in V$, 
where $\epsilon_{1}, \dots, \epsilon_{k}$ 
are arbitrary positive real numbers satisfying 
$\epsilon_{1}>\cdots >\epsilon_{k}$. 
Since $V$ is generated by $X$, elements of the form 
$$e_{k}(u_{1}\otimes \cdots \otimes u_{k}\otimes u\otimes
\mu_{m_{1}, \dots, m_{k}}),$$
$k\ge 0$, 
$u_{1}, \dots, u_{k}, u\in X$ and $m_{1}, \dots, m_{k}\in {\Bbb Z}$,
span $V$. 

Let 
$$h_{1}=e_{k}(u_{1}\otimes \cdots \otimes u_{k}
\otimes u\otimes
\mu_{m_{1}, \dots, m_{k}})$$
and 
$$h_{2}=e_{l}(v_{1}\otimes \cdots \otimes v_{l}
\otimes v\otimes
\mu_{n_{1}, \dots, n_{l}})$$
be two such elements of $V$.
Then 
\begin{eqnarray}\label{3.11}
\lefteqn{\langle \lambda, 
(\overline{\nu}_{Y}([D(z; r_{1}, r_{2})]))(h_{1}\otimes 
h_{2})\rangle}\nno\\
&&=\langle \lambda, e_{k+l+1}(u_{1}\otimes \cdots \otimes 
u_{k}\otimes u\nno\\
&&\hspace{6em} \otimes v_{1}\otimes \cdots \otimes 
v_{l}\otimes  v\otimes \beta_{k,
l}(\mu_{m_{1}, \dots, m_{k}}, 
\mu_{n_{1}, \dots, n_{l}}))\rangle\nno\\
&&=\beta_{k,
l}(\mu_{m_{1}, \dots, m_{k}}, 
\mu_{n_{1}, \dots, n_{l}})\nno\\
&&\hspace{6em}(g_{k+l+1}(
\lambda \otimes u_{1}\otimes \cdots \otimes 
u_{k}\otimes u\otimes v_{1}\otimes \cdots \otimes 
v_{l}\otimes v))\nno\\
&&=\mu_{m_{1}, \dots, m_{k}}
(g_{k}(
\alpha_{l}(\lambda\otimes v_{1}\otimes \cdots \otimes 
v_{l}\otimes   v\otimes \mu_{n_{1}, \dots, n_{l}})
\otimes u_{1}\otimes  \cdots \otimes u_{k}\otimes
u)).\nno\\
&&
\end{eqnarray}
As in the proof of Proposition \ref{3-4},
to avoid notational confusions, we view 
$$g_{k}(
\alpha_{l}(\lambda\otimes v_{1}\otimes \cdots \otimes 
v_{l}\otimes   v\otimes \mu_{n_{1}, \dots, n_{l}})
\otimes u_{1}\otimes  \cdots \otimes u_{k}\otimes
u)$$ 
as an element of $F_{k}^{\zeta}$
and
$$g_{l}(
(P_{n}(Q(u_{1}, \dots, u_{k}, v;
\zeta_{1}, \dots, \zeta_{k}))\diamond_{[D(z; r_{1}, r_{2})]}\lambda)
\otimes v_{1}\otimes \cdots \otimes 
v_{l}\otimes   v)$$
as an element of $F_{l}^{\eta}$. 
We also view $\mu_{m_{1}, \dots, m_{k}}$ and $\nu_{n_{1}, \dots,
n_{l}}$ as elements of $(F_{k}^{\zeta})^{*}$ and $(F_{l}^{\eta})^{*}$,
respectively. 
Then the right-hand side of 
(\ref{3.11}) is equal to 
\begin{eqnarray}\label{3.12}
\lefteqn{\mu_{m_{1}, \dots, m_{k}}
\Biggl(\sum_{n\in {\Bbb Z}}\alpha_{l}(\lambda\otimes 
v_{1}\otimes \cdots \otimes 
v_{l}\otimes   v\otimes \mu_{n_{1}, \dots, n_{l}})}\nno\\
&&\hspace{6em}
(P_{n}(Q(u_{1}, \dots, u_{k}, u;
\zeta_{1}, \dots, \zeta_{k}))))\Biggr)\nno\\
&&=\mu_{m_{1}, \dots, m_{k}}
\Biggl(\sum_{n\in {\Bbb Z}}\mu_{n_{1}, \dots, n_{l}}
(g_{l}( 
(P_{n}(Q(u_{1}, \dots, u_{k}, u;
\zeta_{1}, \dots, \zeta_{k})) \diamond_{[D(z; r_{1}, r_{2})]}\lambda)\nno\\
&&\hspace{6em}
\otimes v_{1}\otimes \cdots \otimes 
v_{l}\otimes   v))\Biggr)\nno\\
&&=\mu_{m_{1}, \dots, m_{k}}
\Biggl(\sum_{n\in {\Bbb Z}}\mu_{n_{1}, \dots, n_{l}}
\Biggl(\sum_{m\in {\Bbb Z}}(P_{n}(Q(u_{1}, \dots, u_{k}, u;
\zeta_{1}, \dots, \zeta_{k})) \diamond_{[D(z; r_{1}, r_{2})]}\lambda)\nno\\
&&\hspace{6em}
(P_{m}(Q(v_{1}, 
\dots, v_{l}, v;
\eta_{1}, \dots, \eta_{l})))\Biggr)\Biggr)\nno\\
&&=\mu_{m_{1}, \dots, m_{k}}
\Biggl(\sum_{n\in {\Bbb Z}}\mu_{n_{1}, \dots, n_{l}}
\Biggl(\sum_{m\in {\Bbb Z}}\sum_{p\in {\Bbb Z}}\lambda(
P_{p}(r_{1}^{L(0)}\cdot \nno\\
&&\hspace{6em}\cdot  Y(P_{n}(Q(u_{1}, 
\dots, u_{k}, u;
\zeta_{1}, \dots, \zeta_{k})), z)\cdot \nno\\
&&\hspace{6em}\cdot r_{2}^{L(0)}
P_{m}(Q(v_{1}, 
\dots, v_{l}, v;
\eta_{1}, \dots, \eta_{l})))\Biggr)\Biggr)\nno\\
&&=\frac{1}{2\pi \sqrt{-1}}\oint_{|\zeta_{1}|=\epsilon_{1}}\cdots 
\frac{1}{2\pi \sqrt{-1}}\oint_{|\zeta_{k}|=\epsilon_{k}}\zeta_{1}^{m_{1}}\cdots
\zeta_{k}^{m_{k}}\cdot \nno\\
&&\hspace{6em}\cdot \sum_{n\in {\Bbb Z}}
\frac{1}{2\pi \sqrt{-1}}\oint_{|\eta_{1}|=\delta_{1}}\cdots 
\frac{1}{2\pi \sqrt{-1}}\oint_{|\eta_{l}|=\delta_{l}}\eta_{1}^{n_{1}}\cdots
\eta_{l}^{n_{l}}\cdot \nno\\
&&\hspace{6em}\cdot \sum_{m\in {\Bbb Z}}\sum_{p\in {\Bbb Z}}\lambda(
P_{p}(r_{1}^{L(0)}Y(P_{n}(Q(u_{1}, 
\dots, u_{k}, u;
\zeta_{1}, \dots, \zeta_{k})), z)\cdot \nno\\
&&\hspace{6em} \cdot r_{2}^{L(0)}
P_{m}(Q(v_{1}, 
\dots, v_{l}, v;
\eta_{1}, \dots, \eta_{l})))\nno\\
&&=\sum_{n\in {\Bbb Z}}\sum_{m\in {\Bbb Z}}
\sum_{p\in {\Bbb Z}}\lambda(P_{p}(r_{1}^{L(0)}Y(P_{n}(\res_{x_{1}}\cdots 
\res_{x_{k}}x_{1}^{m_{1}}\cdots x_{k}^{m_{k}}\cdot \nno\\
&&\hspace{6em}\cdot 
Y(u_{1}, x_{1})
\cdots Y(u_{k}, x_{1})u), z)r_{2}^{L(0)}\cdot \nno\\
&&\hspace{6em}\cdot 
P_{n}(\res_{y_{1}}\cdots 
\res_{y_{l}}y_{1}^{n_{1}}\cdots y_{l}^{n_{l}}
Y(v_{1}, y_{1})
\cdots Y(v_{l}, y_{l})v))\nno\\
&&=\sum_{p\in {\Bbb Z}}\lambda(
P_{p}(Y(r_{1}^{L(0)}h_{1}, z)r_{2}^{L(0)}h_{2}))
\nno\\
&&=\langle \lambda, 
Y(r_{1}^{L(0)}h_{1}, z)r_{2}^{L(0)}h_{2}\rangle.
\end{eqnarray}
From (\ref{3.11}) and (\ref{3.12}), we get 
\begin{eqnarray*}
\lefteqn{(\overline{\nu}_{Y}([D(z; r_{1}, r_{2})]))(h_{1}\otimes 
h_{2})}\nno\\
&&=Y(r_{1}^{L(0)}h_{1}, z)r_{2}^{L(0)}h_{2}\nno\\
&&=(\nu_{Y}((z; {\bf 0}, (1/r_{1}, {\bf 0}), (1/r_{2}, 
{\bf 0}))))(h_{1}\otimes h_{2}),
\end{eqnarray*}
proving the theorem.
\epfv

\renewcommand{\theequation}{\thesection.\arabic{equation}}
\renewcommand{\therema}{\thesection.\arabic{rema}}
\setcounter{equation}{0}
\setcounter{rema}{0}

\section{A locally convex completion of a finitely-gen\-erated module 
and the vertex operator map}

In this section, we discuss modules. Since the constructions and 
proofs are all similar to the case of algebras, we shall only 
state the results and
point out the slight differences between the constructions and proofs 
in the case of modules 
and those in the case of algebras. 

Let $V$ be a ${\Bbb Z}$-graded finitely generated 
vertex algebra satisfying 
the standard grading-restriction axioms and $W$ a ${\Bbb C}$-graded 
finitely-generated $V$-module
satisfying the standard grading-restriction axioms, that is,
$$W=\coprod_{n\in {\Bbb C}}W_{(n)},$$
$$\dim W_{(n)}<\infty,$$
$n\in {\Bbb
C}$ and 
$$W_{(n)}=0$$ for $n\in {\Bbb C}$
whose real part is
sufficiently small.  
Let $M$ be the finite-dimensional vector space
spanned by a set of generators of $W$. 
As in Section 2, we assume that $X$ is a finite-dimensional
subspace containing ${\bf 1}$ of $V$ spanned by a set of generators of $V$. 

We construct spaces $\tilde{G}^{W}$ and $F_{k}^{W}$, 
  $k\ge
0$, and
their
duals  in the same way as the
constructions  of  $\tilde{G}$ and $F_{k}$, 
 $k\ge
0$, and their
duals  in Section
2, except that  $V^{*}$ is replaced by $W^{*}$, 
$v\in V$ is replaced by  $w\in W$ and the vertex
operator map $Y$ for $V$ are replaced by vertex operator map $Y_{W}$
for the module
$W$. We also define linear maps 
$g^{W}_{k}$, $\iota_{F_{k}^{W}}$,  $\gamma_{k}^{W}$, 
$e^{W}_{k}$,
$k\ge 0$, and linear maps induced from them in the same way as 
in the definitions of 
$g_{k}$, $\iota_{F_{k}}$,  $\gamma_{k}$, 
$e_{k}$,
$k\ge 0$, and induced  linear maps in Section 2, except that 
the spaces are replaced by the corresponding spaces involving $W$ and $M$.

For $k\ge 0$,  $e_{k}^{W}$ is a  linear map from $V^{\otimes
k}\otimes W\otimes (F_{k}^{W})^{*}$ to $(\tilde{G}^{W})^{*}$. Let
$$G_{k}^{W}=e_{k}^{W}(X^{\otimes
k}\otimes M\otimes (F_{k}^{W})^{*})$$
 and 
$$G^{W}=
\bigcup_{k\ge 0}G_{k}^{W}.$$
As in the case of $\tilde{G}$ and $H$,
the spaces $\tilde{G}^{W}$ and $(\tilde{G}^{W})^{*}$ form a
dual pair of vector spaces 
and thus we have a locally convex topology on $(\tilde{G}^{W})^{*}$. 
As in Section 2, we have:

\begin{propo}
For any $k\ge 0$, the linear map 
$$e^{W}_{k}|_{X^{\otimes
k}\otimes M\otimes (F_{k}^{W})^{*}}: X^{\otimes
k}\otimes M\otimes (F_{k}^{W})^{*}\to (\tilde{G}^{W})^{*}$$ is continuous. In
particular, 
$$X^{\otimes
k}\otimes M\otimes (F_{k}^{W})^{*}/(e^{W}_{k}|_{X^{\otimes
k}\otimes M\otimes (F_{k}^{W})^{*}})^{-1}(0)$$
and  $G^{W}_{k}$ is a 
locally convex space.\epf
\end{propo}

For $k\ge 0$, let $H_{k}^{W}$ be the completion of $G_{k}^{W}$ and 
$$H^{W}=\bigcup_{k\ge 0}H_{k}^{W}.$$
Then
$H_{k}^{W}$, $k\ge 0$,  are complete locally convex space. 

\begin{propo}
For $k\ge 0$, $H^{W}_{k}$ can be embedded canonically 
in $H^{W}_{k+1}$ and the topology on
$H^{W}_{k}$ is induced from the topology on $H^{W}_{k+1}$. \epf
\end{propo}

\begin{theo}
The vector space $H^{W}$ equipped with the strict inductive limit
topology is a locally convex completion of $W$. In particular, $G^{W}$
is dense in $H^{W}$.\epf
\end{theo}

Next we extend the map 
$$Y_{W}(r_{1}^{L(0)}\cdot, z)r_{2}^{L(0)}\cdot : 
V\otimes W
\to \overline{W}$$ 
associated to $[D(z;  r_{1},
r_{2})]$
to 
a linear map $\overline{\nu}_{Y_{W}}([D(z;  r_{1},
r_{2})])$ from $H\widetilde{\otimes}H^{W}$ to $H^{W}$.

We define 
$u\diamond_{[D(z; r_{1}, r_{2})]}\lambda\in W^{*}$ 
for $\lambda\in \tilde{G}^{W}$ and $u\in V$,  the maps
$$\alpha^{W}_{l}: \tilde{G}^{W}\otimes X^{l}\otimes W\otimes 
(F_{l}^{W})^{*}\to \tilde{G},$$
$l\ge 0$, and 
$$\beta^{W}_{k, l}: F^{*}_{k}\otimes (F_{l}^{W})^{*}
\to (F^{W}_{k+l+1})^{*},$$
$k, l\ge 0$, in the same way
as in the definitions in Section 3 of 
$u\diamond_{[D(z; r_{1}, r_{2})]}\lambda\in V^{*}$ for 
$\lambda\in \tilde{G}$ and $u\in V$, the maps
$$\alpha_{l}: \tilde{G}^{W}\otimes  X^{l+1}\otimes 
F^{*}_{l}\to V^{*},$$
$l\ge 0$, and 
$$\beta_{k, l}: F^{*}_{k}\otimes F^{*}_{l}
\to F^{*}_{k+l+1},$$
$k, l\ge 0$. (Note that the image of $\alpha^{W}_{l}$ is 
in $\tilde{G}$,  not in $\tilde{G}^{W}$. This is the reason why the domain
of 
$\beta^{W}_{k, l}$ is $F^{*}_{k}\otimes (F^{W}_{l})^{*}$,  not 
$(F^{W}_{k})^{*}\otimes (F^{W}_{l})^{*}$.)

Let 
$$h_{1}=e_{k}(u_{1}\otimes \cdots \otimes u_{k}\otimes
u\otimes \mu)\in G_{k}$$
and 
$$h_{2}=e^{W}_{l}(v_{1}\otimes \cdots \otimes v_{l}
\otimes
w\otimes \nu)\in G^{W}_{l}$$
where $u_{1}, \dots, u_{k}, u, v_{1}, \dots, 
v_{l}\in X$, $w\in M$ and 
$\mu\in F^{*}_{k}$, $\nu\in 
(F^{W}_{l})^{*}$. Define 
\begin{eqnarray*}
\lefteqn{(\overline{\nu}_{Y_{W}}([D(z; r_{0}, r_{1},
r_{2})]))(h_{1}\otimes h_{2})}\nno\\
&&=e^{W}_{k+l+1}(u_{1}\otimes \cdots \otimes 
u_{k}\otimes u\otimes v_{1}\otimes \cdots \otimes 
v_{l}\otimes  w\otimes \beta^{W}_{k, l}(\mu, \nu)).
\end{eqnarray*}
Note that any element of $G_{k}$ or $G^{W}_{l}$ is a linear
combination of elements of the form
$h_{1}$ or $h_{2}$, respectively, above, and that $k$ and $l$ are
arbitrary. Thus we obtain a linear map 
$$\overline{\nu}_{Y_{W}}([D(z; r_{0}, r_{1},
r_{2})])\mbar_{G\otimes G^{W}}: G\otimes G^{W}\to G^{W}.$$

\begin{propo}
The map $\overline{\nu}_{Y_{W}}([D(z; r_{1},
r_{2})])\mbar_{G\otimes G^{W}}$ is continuous.\epf
\end{propo}

By this proposition and the fact that $G$ and $G^{W}$ are dense in $H$
and $H^{W}$, respectivel, 
we can extend $\overline{\nu}_{Y_{W}}([D(z; r_{1},
r_{2})])\mbar_{G\otimes G^{W}}$ to a continuous linear map
$\overline{\nu}_{Y_{W}}([D(z; r_{1},
r_{2})])$
from $H\widetilde{\otimes} H^{W}$ to $H^{W}$. 

\begin{theo}
The map $\overline{\nu}_{Y_{W}}([D(z; r_{1},
r_{2})])$ is a continuous extension of 
$$Y_{W}(r_{1}^{L(0)}\cdot, z)r_{2}^{L(0)}\cdot : V\otimes W\to \overline{W}$$
to $H\widetilde{\otimes}H^{W}$. \epf
\end{theo}

{\small \sc Department of Mathematics, Rutgers University,
110 Frelinghuysen Rd., Piscataway, NJ 08854-8019}

{\em E-mail address}: yzhuang@math.rutgers.edu

\end{document}